\documentclass[letterpaper]{IEEEtran}
\usepackage{makeidx}
\usepackage{graphicx}
\usepackage{multicol}
\usepackage[bottom]{footmisc}
\usepackage{subfigure}
\usepackage{cite,color,comment,xspace}
\usepackage{amsfonts}
\usepackage{amsmath}
\usepackage{mathrsfs}
\usepackage{amssymb}
\usepackage{algorithm}
\usepackage{algorithmic}
\usepackage{url}
\usepackage{times}
\usepackage{mathptmx}
\usepackage{enumerate}
\usepackage{epstopdf}
\usepackage{subfigure}
\usepackage{epsfig}%
\newtheorem{theorem}{Theorem}

\newtheorem{proposition}[theorem]{Proposition}

\makeatletter

\@addtoreset{theorem}{section}

\@addtoreset{corollary}{section}

\@addtoreset{lemma}{section}

\@addtoreset{definition}{section}
\makeatother
\begin{document}

\title{{\LARGE \textbf{An Optimal Control Approach to the Multi-Agent Persistent
Monitoring Problem in Two-Dimensional Spaces}}}
\author{\textbf{Xuchao Lin} and \textbf{Christos G. Cassandras}
\thanks{{\footnotesize The authors' work is supported in part by NSF under
Grant CNS-1239021, by AFOSR under grant FA9550-12-1-0113, by ONR under grant
N00014-09-1-1051, and by ARO under Grant W911NF-11-1-0227.}} \thanks{*
Division of Systems Engineering and Center for Information and Systems
Engineering, Boston University; e-mail:\emph{\{mmxclin,cgc\}@bu.edu}.}}
\maketitle

\begin{abstract}
We address the persistent monitoring problem in two-dimensional mission spaces
where the objective is to control the trajectories of multiple cooperating
agents to minimize an uncertainty metric. In a one-dimensional mission space,
we have shown that the optimal solution is for each agent to move at maximal
speed and switch direction at specific points, possibly waiting some time at
each such point before switching. In a two-dimensional mission space, such
simple solutions can no longer be derived. An alternative is to optimally
assign each agent a linear trajectory, motivated by the one-dimensional
analysis. We prove, however, that elliptical trajectories outperform linear
ones. With this motivation, we formulate a parametric optimization problem in
which we seek to determine such trajectories. We show that the problem can be
solved using Infinitesimal Perturbation Analysis (IPA) to obtain performance
gradients on line and obtain a complete and scalable solution. Since the
solutions obtained are generally locally optimal, we incorporate a stochastic
comparison algorithm for deriving globally optimal elliptical trajectories.
Numerical examples are included to illustrate the main result, allow for
uncertainties modeled as stochastic processes, and compare our proposed
scalable approach to trajectories obtained through off-line computationally
intensive solutions.

\end{abstract}

\section{Introduction}

Autonomous cooperating agents may be used to perform tasks such as coverage
control \cite{zhong2011distributed,cortes2004coverage}, surveillance
\cite{grocholsky2006cooperative} and environmental sampling
\cite{smith2011persistent,paley2008cooperative,Dames2012adecentralized}.
\emph{Persistent monitoring} (also called \textquotedblleft persistent
surveillance\textquotedblright\ or \textquotedblleft persistent
search\textquotedblright) arises in a large dynamically changing environment
which cannot be fully covered by a stationary team of available agents. Thus,
persistent monitoring differs from traditional coverage tasks due to the
perpetual need to cover a changing environment, i.e., all areas of the mission
space must be sensed infinitely often. The main challenge in designing control
strategies in this case is in balancing the presence of agents in the changing
environment so that it is covered over time optimally (in some well-defined
sense) while still satisfying sensing and motion constraints.

Control and motion planning for agents performing persistent monitoring tasks
have been studied in the literature, e.g., see
\cite{SmithTRO2010Submit,Soltero2011collision,nigam2008persistent,hokayem2008persistent,Julian2012nonparametric,chen2012multi,Lan2013planning}%
. In \cite{cassandras2013optimal}, we addressed the persistent monitoring
problem by proposing an \emph{optimal control} framework to drive multiple
cooperating agents so as to minimize a metric of uncertainty over the
environment. This metric is a function of both space and time such that
uncertainty at a point grows if it is not covered by any agent sensors. To
model sensor coverage, we define a probability of detecting events at each
point of the mission space by agent sensors. Thus, the uncertainty of the
environment decreases with a rate proportional to the event detection
probability, i.e.,\textit{ } the higher the sensing effectiveness is, the
faster the uncertainty is reduced. It was shown in
\cite{cassandras2013optimal} that the optimal control problem can be reduced
to a \emph{parametric} optimization problem. In particular, the optimal
trajectory of each agent is to move at full speed until it reaches some
switching point, dwell on the switching point for some time (possibly zero),
and then switch directions. Thus, each agent's optimal trajectory is fully
described by a set of switching points $\{\theta_{1},\ldots,\theta_{K}\}$ and
associated waiting times at these points, $\{w_{1},\ldots,w_{K}\}$. This
allows us to make use of Infinitesimal Perturbation Analysis (IPA)
\cite{cassandras2009perturbation}\textbf{ }to determine gradients of the
objective function with respect to these parameters and subsequently obtain
optimal switching locations and waiting times that fully characterize an
optimal solution. It also allows us to exploit robustness properties of IPA to
readily extend this solution approach to a \emph{stochastic} uncertainty model.

In this paper, we address the same persistent monitoring problem in a
two-dimensional (2D) mission space. Using an analysis similar to the
one-dimensional (1D) case, we find that we can no longer identify a parametric
representation of optimal agent trajectories. A complete solution requires a
computationally intensive process for solving a Two Point Boundary Value
Problem (TPBVP) making any on-line solution to the problem infeasible.
Motivated by the simple structure of the 1D problem, it has been suggested to
assign each agent a linear trajectory for which the explicit 1D solution can
be used. One could then reduce the problem to optimally carrying out this
assignment. However, in a 2D space it is not obvious that a linear trajectory
is a desirable choice. Indeed, a key contribution of this paper is to formally
prove that an elliptical agent trajectory outperforms a linear one in terms of
the uncertainty metric we are using. Motivated by this result, we formulate a
2D persistent monitoring problem as one of determining optimal elliptical
trajectories for a given number of agents, noting that this includes the
possibility that one or more agents share the same trajectory. We show that
this problem can be explicitly solved using similar IPA techniques as in our
1D analysis. In particular, we use IPA to determine on line the gradient of
the objective function with respect to the parameters that fully define each
elliptical trajectory (center, orientation and length of the minor and major
axes). This approach is scalable in the number of observed events, not states,
of the underlying hybrid system characterizing the persistent monitoring
process, so that it is suitable for on-line implementation. However, the
standard gradient-based optimization process we use is generally limited to
local, rather than global optimal solutions. Thus, we adopt a stochastic
comparison algorithm from the literature \cite{bao1996stochastic} to overcome
this problem.

Section \ref{sec:problemformulation} formulates the optimal control problem
for 2D mission spaces and Section \ref{sec:optimalsolution} presents the
solution approach. In Section IV we establish our key result that elliptical
agent trajectories outperform linear ones in terms of minimizing an
uncertainty metric per unit area. In Section V we formulate and solve the
problem of determining optimal elliptical agent trajectories using an
algorithm driven by gradients evaluated through IPA. In Section VI we
incorporate a stochastic comparison algorithm for obtaining globally optimal
solutions and in Section VII we provide numerical results to illustrate our
approach and compare it to computationally intensive solutions based on a
TPBVP solver. Section \ref{sec:concl} concludes the paper.

\section{Persistent Monitoring Problem Formulation}

\label{sec:problemformulation} We consider $N$ mobile agents in a 2D
rectangular mission space $\Omega\equiv\lbrack0,L_{1}]\times\lbrack
0,L_{2}]\subset\mathbb{R}^{2}$. Let the position of the agents at time $t$ be
$s_{n}(t)=[s_{n}^{x}(t),s_{n}^{y}(t)]$ with $s_{n}^{x}(t)\in\lbrack0,L_{1}]$
and $s_{n}^{y}(t)\in\lbrack0,L_{2}]$, $n=1,\ldots,N$, following the dynamics:%
\begin{equation}
\dot{s}_{n}^{x}(t)=u_{n}(t)\cos\theta_{n}\left(  t\right)  ,\text{ \ \ }%
\dot{s}_{n}^{y}(t)=u_{n}\left(  t\right)  \sin\theta_{n}\left(  t\right)
\label{eq:DynS}%
\end{equation}
where $u_{n}\left(  t\right)  $ is the scalar speed of the $n$th agent and
$\theta_{n}\left(  t\right)  $ is the angle relative to the positive direction
that satisfies $0\leq\theta_{n}\left(  t\right)  <2\pi$. Thus, we assume that
each agent controls its orientation and speed. Without loss of generality,
after some rescaling of the size of the mission space, we further assume that
the speed is constrained by $0\leq u_{n}(t)\leq1$, $n=1,\ldots,N$. Each agent
is represented as a particle in the 2D space, thus we ignore the case of two
or more agents colliding with each other.

We associate with every point $[x,y]\in\Omega$ a function $p_{n}(x,y,s_{n})$
that measures the probability that an event at location $[x,y]$ is detected by
agent $n$. We also assume that $p_{n}(x,y,s_{n})=1$ if $[x,y]=s_{n}$, and that
$p_{n}(x,y,s_{n})$ is monotonically nonincreasing in the Euclidean distance
$D(x,y,s_{n})\equiv||[x,y]-s_{n}||$ between $[x,y]$ and $s_{n}$, thus
capturing the reduced effectiveness of a sensor over its range which we
consider to be finite and denoted by $r_{n}$ (this is the same as the concept
of \textquotedblleft sensor footprint\textquotedblright\ commonly used in the
robotics literature.) Therefore, we set $p_{n}(x,y,s_{n})=0$ when
$D(x,y,s_{n})>r_{n}$. Our analysis is not affected by the precise sensing
model $p_{n}(x,y,s_{n})$, but we mention here as an example the linear decay
model used in \cite{cassandras2013optimal}:%
\begin{equation}
p_{n}(x,y,s_{n})=\left\{
\begin{array}
[c]{cc}%
\frac{1}{C}(1-\frac{D(x,y,s_{n})}{r_{n}}), & \text{if }D(x,y,s_{n})\text{ }\leq r_{n}\\
0, & \text{if }D(x,y,s_{n})\text{ }>r_{n}%
\end{array}
\right.  \label{eq:multiLinearModel}%
\end{equation}
where $C$ is a normalization constant. Next, consider a set of points $\{[\alpha_{i},\beta_{i}]$, $i=1,\ldots,M\}$,
$[\alpha_{i},\beta_{i}]\in\Omega$, and associate a time-varying measure of
uncertainty with each point $[\alpha_{i},\beta_{i}]$, which we denote by
$R_{i}(t)$.
The set of points $\{[\alpha_{1},\beta_{1}],\ldots,[\alpha_{M},\beta_{M}]\}$
may be selected to contain specific \textquotedblleft points of
interest\textquotedblright\ in the environment, or simply to sample points in
the mission space. Alternatively, we may consider a partition of $\Omega$ into
$M$ rectangles denoted by $\Omega_{i}$ whose center points are $[\alpha
_{i},\beta_{i}]$. We can then set $p_{n}(x,y,s_{n}(t))=p_{n}(\alpha_{i}%
,\beta_{i},s_{n}\left(  t\right)  )$ for all $\{[x,y]|[x,y]\in\Omega
_{i},[\alpha_{i},\beta_{i}]\in\Omega_{i}\}$, i.e., for all $[x,y]$ in the
rectangle $\Omega_{i}$ that $[\alpha_{i},\beta_{i}]$ belongs to. 
In order to avoid the uninteresting case where there is a large
number of agents who can adequately cover the mission space, we assume that for
any $\mathbf{s}(t)$, there
exists some point $[x,y]\in\Omega$ with $P(x,y,\mathbf{s}(t))=0.$ This means
that for any assignment of $N$ agents at time $t$, there is always at least
one point in the mission space that cannot be sensed by any agent.
Therefore,
the joint probability of detecting an event at location $[\alpha_{i},\beta
_{i}]$ by all the $N$ agents (assuming detection independence) is%
\[
P_{i}\left(  \mathbf{s}(t)\right)  =1-{\textstyle\prod\limits_{n=1}^{N}%
}\left[  1-p_{n}(\alpha_{i},\beta_{i},s_{n}\left(  t\right)  )\right]
\]
where we set $\mathbf{s}(t)=[s_{1}\left(  t\right)  ,\ldots,s_{N}\left(
t\right)  ]^{\text{T}}$. Similar to the 1D analysis in
\cite{cassandras2013optimal}, we define uncertainty functions $R_{i}(t)$
associated with the rectangles $\Omega_{i}$, $i=1,\ldots,M$, so that they have
the following properties: $(i)$ $R_{i}(t)$ increases with a prespecified rate
$A_{i}$ if $P_{i}\left(  \mathbf{s}(t)\right)  =0$, $(ii)$ $R_{i}(t)$
decreases with a fixed rate $B - A_{i}$ if $P_{i}\left(  \mathbf{s}(t)\right)  =1$ and
$(iii)$ $R_{i}(t)\geq0$ for all $t$. It is then natural to model uncertainty
so that its decrease is proportional to the probability of detection. In
particular, we model the dynamics of $R_{i}(t)$, $i=1,\ldots,M$, as follows:
\begin{equation}
\dot{R}_{i}(t)=\left\{
\begin{array}
[c]{ll}%
0 & \text{if }R_{i}(t)=0,\text{ }A_{i}\leq BP_{i}\left(  \mathbf{s}(t)\right)
\\
A_{i}-BP_{i}\left(  \mathbf{s}(t)\right)  & \text{otherwise}%
\end{array}
\right.  \label{eq:multiDynR}%
\end{equation}
where we assume that initial conditions $R_{i}(0)$, $i=1,\ldots,M$, are given
and that $B>A_{i}>0$ for all $i=1,\ldots,M$; thus, the uncertainty strictly
decreases when there is perfect sensing $P_{i}\left(  \mathbf{s}(t)\right)
=1$.

The goal of the optimal persistent monitoring problem we consider is to
control through $u_{n}\left(  t\right)  $, $\theta_{n}(t)$ in (\ref{eq:DynS})
the movement of the $N$ agents so that the cumulative uncertainty over all
sensing points $\{[\alpha_{1},\beta_{1}],\ldots,[\alpha_{M},\beta_{M}]\}$ is
minimized over a fixed time horizon $T$. Thus, setting $\mathbf{u}\left(
t\right)  =\left[  u_{1}\left(  t\right)  ,\ldots,u_{N}\left(  t\right)
\right]  $ and $\mathbf{\theta}(t)=\left[  \theta_{1}\left(  t\right)
,\ldots,\theta_{N}\left(  t\right)  \right]  $ we aim to solve the following
optimal control problem\textbf{ P1}:
\begin{equation}
\mathbf{P1:}\text{ \ \ \ \ \ \ \ }\min_{\mathbf{u}\left(  t\right)
,\mathbf{\theta}(t)}\text{ \ }J=\int_{0}^{T}\sum_{i=1}^{M}R_{i}(t)dt
\label{eq:costfunction}%
\end{equation}
subject to the agent dynamics \eqref{eq:DynS}, uncertainty dynamics
\eqref{eq:multiDynR}, control constraint $0\leq u_{n}(t)\leq1$, $0\leq
\theta_{n}(t)\leq2\pi$, $t\in\lbrack0,T]$, and state constraints $s_{n}%
(t)\in\Omega$ for all $t\in\lbrack0,T]$, $n=1,\ldots,N$.

\textbf{Remark 1}. The modeling of the uncertainty value $R_{i}(t)$ in a 2D
environment is a direct extension of \cite{cassandras2013optimal} in the 1D
environment setting where it was described how persistent monitoring can be
viewed as a polling system, with each rectangle $\Omega_{i}$ associated with a
\textquotedblleft virtual queue\textquotedblright\ where uncertainty
accumulates with inflow rate $A_{i}$. Each agent acts as a \textquotedblleft
server\textquotedblright\ visiting these virtual queues with a time-varying
service rate given by $BP_{i}\left(  \mathbf{s}(t)\right)  $, controllable
through all agent positions at time $t$. Metrics other than
(\ref{eq:costfunction}) are of course possible, e.g., maximizing the mutual
information or minimizing the spectral radius of the error covariance matrix
\cite{Zhang2010on} if specific \textquotedblleft point of
interest\textquotedblright\ are identified with known properties.

\section{Optimal Control Solution}

\label{sec:optimalsolution} We first characterize the optimal control solution
of problem\textbf{ P1}. We define the state vector $\mathbf{x}\left(
t\right)  =[s_{1}^{x}\left(  t\right)  ,s_{1}^{y}\left(  t\right)
,\ldots,s_{N}^{x}(t),s_{N}^{y}(t),R_{1}\left(  t\right)  ,\ldots,R_{M}\left(
t\right)  ]^{\mathtt{T}}$ and the associated costate vector $\mathbf{\lambda
}\left(  t\right)  =[\mu_{1}^{x}(t),\mu_{1}^{y}(t),\ldots,\mu_{N}^{x}%
(t),\mu_{N}^{y}(t),\lambda_{1}\left(  t\right)  ,\ldots,\lambda_{M}\left(
t\right)  ]^{\mathtt{T}}$. In view of the discontinuity in the dynamics of
$R_{i}(t)$ in (\ref{eq:multiDynR}), the optimal state trajectory may contain a
boundary arc when $R_{i}(t)=0$ for any $i$; otherwise, the state evolves in an
interior arc \cite{bryson1975applied}. 
This follows from the fact, proved in 
 \cite{cassandras2013optimal} and
\cite{lin2013optimal} that it is never optimal for agents to reach the mission
space boundary.
We analyze the system operating
in such an interior arc and omit the state constraint $s_{n}(t)\in\Omega,$
$n=1,\ldots,N,$ $t\in\lbrack0,T]$. Using (\ref{eq:DynS}) and
(\ref{eq:multiDynR}), the Hamiltonian is%
\begin{align}
H &  =\sum_{i}R_{i}(t)+\sum_{i}\lambda_{i}\dot{R}_{i}(t)\nonumber\\
&  +\sum_{n}\mu_{n}^{x}\left(  t\right)  u_{n}\left(  t\right)  \cos\theta
_{n}\left(  t\right)  +\sum_{n}\mu_{n}^{y}\left(  t\right)  u_{n}\left(
t\right)  \sin\theta_{n}\left(  t\right)  \label{Hamiltonian}%
\end{align}
and the costate equations $\mathbf{\dot{\lambda}}=-\frac{\partial H}{\partial \mathbf{x}}$ are%
\begin{align}
\dot{\lambda}_{i}(t)& =-\frac{\partial H}{\partial R_{i}}%
=-1\label{eq:multiDynCoI}\\%
\dot{\mu}_{n}^{x}\left(  t\right)   &  =-\frac{\partial H}{\partial s_{n}^{x}%
}=-\sum_{i}\frac{\partial}{\partial s_{n}^{x}}\lambda_{i}\dot{R}%
_{i}(t)\nonumber\\
&  =-\sum_{[\alpha_{i},\beta_{i}]\in\mathcal{R}\left(  s_{n}\right)  }%
\frac{B\lambda_{i}\left(  s_{n}^{x}-\alpha_{i}\right)  }{r_{n}D(\alpha
_{i},\beta_{i},s_{n}(t))}\textstyle\prod\limits_{d\neq n}^{N}\left[
1-p_{d}\left(  \omega_{i},s_{d}\left(  t\right)  \right)  \right]  \label{eq:DynOfLamX}\\
\dot{\mu}_{n}^{y}\left(  t\right)   &  =-\frac{\partial H}{\partial s_{n}^{y}%
}=-\sum_{i}\frac{\partial}{\partial s_{n}^{y}}\lambda_{i}\dot{R}%
_{i}(t)\nonumber\\
&  =-\sum_{[\alpha_{i},\beta_{i}]\in\mathcal{R}\left(  s_{n}\right)  }%
\frac{B\lambda_{i}\left(  s_{n}^{y}-\beta_{i}\right)  }{r_{n}D(\alpha
_{i},\beta_{i},s_{n}(t))}\textstyle\prod\limits_{d\neq n}^{N}\left[
1-p_{d}\left(  \omega_{i},s_{d}\left(  t\right)  \right)  \right]   \label{eq:DynOfLamY}\\
& \nonumber
\end{align}
where $\mathcal{R}\left(  s_{n}\right)  \equiv\left\{  \lbrack\alpha_{i}%
,\beta_{i}]\left\vert D(\alpha_{i},\beta_{i},s_{n})\leq r_{n},\text{ }%
R_{i}>0\right.  \right\}  $ identifies all points $[\alpha_{i},\beta_{i}]$
within the sensing range of the agent using the model in
(\ref{eq:multiLinearModel}). Since we impose no terminal state constraints,
the boundary conditions are $\lambda_{i}(T)=0,$ $i=1,\ldots,M$ and $\mu
_{n}^{x}(T)=0,$ $\mu_{n}^{y}(T)=0,$ $n=1,\ldots,N$.
The implication of (\ref{eq:multiDynCoI}) with $\lambda_{i}\left(  T\right)
=0$ is that $\lambda_{i}\left(  t\right)  =T-t$ for all $t\in\lbrack0,T]$,
$i=1,\ldots,M$ and that $\lambda_{i}\left(  t\right)  $ is monotonically
decreasing starting with $\lambda_{i}\left(  0\right)  =T$. However, this is
only true if the entire optimal trajectory is an interior arc, i.e., all
$R_{i}(t)\geq0$ constraints for all $i=1,\ldots,M$ remain inactive.
We have shown in \cite{cassandras2013optimal} that $\lambda_{i}\left(
t\right)  \geq0$, $i=1,\ldots,M$, $t\in\left[  0,T\right]  $ with equality
holding only if $t=T,$ or $t=$ $t_{0}^{-}$ with $R_{i}\left(  t_{0}\right)
=0$, $R_{i}\left(  t^{\prime}\right)  >0$, where $t^{\prime}\in\lbrack
t_{0}-\delta,t_{0})$, $\delta>0.$ Although this argument holds for the 1D
problem formulation, the proof can be\ directly extended to this 2D
environment. However, the actual evaluation of the full costate vector over
the interval $[0,T]$ requires solving (\ref{eq:DynOfLamX}) and
(\ref{eq:DynOfLamY}), which in turn involves the determination of all points
where the state variables $R_{i}(t)$ reach their minimum feasible values
$R_{i}(t)=0$, $i=1,\ldots,M$. This generally involves the solution of a TPBVP.

From (\ref{Hamiltonian}), after some algebraic operations, we get
\begin{align}
H &  =\sum_{i}R_{i}(t)+\sum_{i}\lambda_{i}\dot{R}_{i}(t)\nonumber \\
& +\sum_{n}%
u_{n}(t)\left[  \mu_{n}^{x}\left(  t\right)  \cos\theta_{n}\left(  t\right)
+\mu_{n}^{y}\left(  t\right)  \sin\theta_{n}\left(  t\right)  \right]
\nonumber\\
&  =\sum_{i}R_{i}(t)+\sum_{i}\lambda_{i}\dot{R}_{i}(t)+\sum_{n}\text{sgn}%
(\mu_{n}^{y}(t))\sqrt{(\mu_{n}^{x}\left(  t\right)  )^{2}+(\mu_{n}%
^{y}\left(  t\right)  )^{2}}\nonumber\\
&  \times u_{n}(t)\left[  \frac{\text{sgn}(\mu_{n}^{y}(t))\mu_{n}^{x}\left(  t\right)\cos\theta_{n}\left(  t\right)
}{\sqrt{(\mu_{n}^{x}\left(  t\right)  )^{2}+(\mu_{n}^{y}\left(  t\right)
)^{2}}} +\frac{|\mu_{n}^{y}\left(  t\right)
|\sin\theta_{n}\left(  t\right)}{\sqrt{(\mu_{n}^{x}\left(  t\right)  )^{2}+(\mu_{n}^{y}\left(  t\right)
)^{2}}}  \right]
\end{align}
where sgn$(\cdot)$ is the sign function. Combining the trigonometric function
terms, we obtain
\begin{align}
H & =\sum_{i}R_{i}(t)+\sum_{i}\lambda_{i}\dot{R}_{i}(t) \nonumber \\
& +\sum_{n}\text{sgn}%
(\mu_{n}^{y}(t))u_{n}\left(  t\right)  \sqrt{(\mu_{n}^{x}\left(  t\right)
)^{2}+(\mu_{n}^{y}\left(  t\right)  )^{2}}\sin(\theta_{n}\left(  t\right)
+\psi_{n}\left(  t\right)  )\label{eq:twoDHamiltonian}%
\end{align}
and $\psi_{n}(t)$ is defined so that $\tan\psi_{n}\left(  t\right)
=\frac{\mu_{n}^{x}\left(  t\right)  }{\mu_{n}^{y}\left(  t\right)  }$ for
$\mu_{n}^{y}(t)\neq0$ and
\[
\psi_{n}(t)=\left\{
\begin{array}
[c]{cc}%
-\frac{\pi}{2}, & \text{if }\mu_{n}^{x}\left(  t\right)  <0\\
\frac{\pi}{2}, & \text{if }\mu_{n}^{x}\left(  t\right)  >0
\end{array}
\right.
\]
for $\mu_{n}^{y}(t)=0$. In what follows, we exclude the case where $\mu_{n}%
^{x}(t)=0$ and $\mu_{n}^{y}(t)=0$ at the same time for any given $n$ over any
finite \textquotedblleft singular interval.\textquotedblright\ Applying the
Pontryagin minimum principle to (\ref{eq:twoDHamiltonian}) with $u_{n}^{\ast
}(t),$ $\theta_{n}^{\ast}(t)$, $t\in\lbrack0,T)$, denoting optimal
controls, we have%
\[
H\left(  \mathbf{x}^{\ast},\mathbf{\lambda}^{\ast},\mathbf{u}^{\ast
},\mathbf{\theta}^{\ast}\right)  =\min_{\mathbf{u}\in\lbrack
0,1]^{N},\mathbf{\theta}\in\lbrack0,2\pi]^{N}}H\left(  \mathbf{x},\mathbf{\lambda
},\mathbf{u},\theta\right)
\]
and it is immediately obvious that it is necessary for an optimal control to
satisfy:%
\begin{equation}
u_{n}^{\ast}\left(  t\right)  =1\label{optu}%
\end{equation}
and
\begin{equation}
\left\{
\begin{array}
[c]{c}%
\sin\left(  \theta_{n}^{\ast}\left(  t\right)  +\psi_{n}\left(  t\right)
\right)  =1\text{, }\ \ \text{if }\mu_{n}^{y}(t)<0\\
\sin\left(  \theta_{n}^{\ast}\left(  t\right)  +\psi_{n}\left(  t\right)
\right)  =-1\text{, \ if }\mu_{n}^{y}(t)>0
\end{array}
\right.  \label{eq:twoDoptimalControl}%
\end{equation}
Note $u_{n}\left(  t\right)  =0$ is not an optimal solution, since we can
always set control $\theta_{n}\left(  t\right)  $ to enforce $\text{sgn}(\mu_{n}^{y}(t)) \sin\left(
\theta_{n}\left(  t\right)  +\psi_{n}\left(  t\right)  \right)  <0$. Thus,
we have%
\begin{equation}
\left\{
\begin{array}
[c]{c}%
\theta_{n}^{\ast}\left(  t\right)  =\frac{\pi}{2}-\psi_{n}\left(  t\right)
\text{, \ \ if }\mu_{n}^{y}(t)<0\\
\theta_{n}^{\ast}\left(  t\right)  =\frac{3\pi}{2}-\psi_{n}\left(
t\right)  \text{, \ \ if }\mu_{n}^{y}(t)>0
\end{array}
\right.  \label{OptTheta1}%
\end{equation}
Clearly, when $\mu_{n}^{y}(t)<0$, the $n$th agent heading is $\theta_{n}%
^{\ast}\left(  t\right)  =\frac{1}{2}\pi-\psi_{n}\left(  t\right)
\in(0,\pi)$ and the agent will move upward in $\Omega$; similarly,
when $\mu_{n}^{y}(t)>0$ the agent will move downward. When $\mu
_{n}^{y}(t)=0$, we have
\[
\psi_{n}(t)=\left\{
\begin{array}
[c]{cc}%
-\frac{\pi}{2}, & \text{if }\mu_{n}^{x}\left(  t\right)  <0\\
\frac{\pi}{2}, & \text{if }\mu_{n}^{x}\left(  t\right)  >0
\end{array}
\right.  \text{ \ and \ }\theta_{n}^{\ast}(t)=\left\{
\begin{array}
[c]{cc}%
0, & \text{if }\mu_{n}^{x}\left(  t\right)  <0\\
\pi, & \text{if }\mu_{n}^{x}\left(  t\right)  >0
\end{array}
\right.
\]
so that the $n$th agent will move horizontally. By symmetry, the agent will
move towards the right when $\mu_{n}^{x}(t)<0$, towards the left when $\mu
_{n}^{x}(t)>0$, and vertically when $\mu_{n}^{x}(t)=0.$ Note that this is
analogous to the 1D problem in \cite{cassandras2013optimal} where the costate
$\lambda_{s_{n}}(t)<0$ implies $u_{n}(t)=1$ and $\lambda_{s_{n}}(t)>0$ implies
$u_{n}(t)=-1$.

Returning to the Hamiltonian in (\ref{Hamiltonian}), the optimal heading
$\theta_{n}^{\ast}(t)$ can be obtained by requiring $\frac{\partial H^{\ast}%
}{\partial\theta_{n}^{\ast}}=0$:%
\[
\frac{\partial H}{\partial\theta_{n}}=-\mu_{n}^{x}(t)u_{n}\left(  t\right)
\sin\theta_{n}\left(  t\right)  +\mu_{n}^{y}(t)u_{n}\left(  t\right)  \cos
\theta_{n}\left(  t\right)  =0
\]
which gives:
\begin{equation}
\tan\theta_{n}^{\ast}(t)=\frac{\mu_{n}^{y}(t)}{\mu_{n}^{x}(t)}
\label{OptTheta2}%
\end{equation}
Applying the tangent operation to both sides of (\ref{OptTheta1}), we can see
that (\ref{OptTheta1}) and (\ref{OptTheta2}) are equivalent to each other.

Since we have shown that $u_{n}^{\ast}\left(  t\right)  =1,$ $n=1,\ldots,N$ in
(\ref{OptTheta1}), we are only left with the task of determining $\theta
_{n}^{\ast}(t),$ $n=1,\ldots,N$. This can be accomplished by solving a
standard TPBVP involving forward and backward integrations of the
state and costate equations to evaluate $\frac{\partial H}{\partial\theta_{n}%
}$ after each such iteration and using a gradient descent approach until the
objective function converges to a (local) minimum. Clearly, this is a
computationally intensive process which scales poorly with the number of
agents and the size of the mission space. In addition, it requires
discretizing the mission time $T$ and calculating every control at each time
step which adds to the computational complexity.

\section{Linear vs Elliptical Agent Trajectories}

Given the complexity of the TPBVP required to obtain an optimal solution of
problem\textbf{ P1}, we seek alternative approaches which may be suboptimal
but are tractable and scalable. The first such effort is motivated by the
results obtained in our 1D analysis, where we found that on a mission space
defined by a line segment $[0,L]$ the optimal trajectory for each agent is to
move at full speed until it reaches some switching point, dwell on the
switching point for some time (possibly zero), and then switch directions.
Thus, each agent's optimal trajectory is fully described by a set of switching
points $\{\theta_{1},\ldots,\theta_{K}\}$ and associated waiting times at
these points, $\{w_{1},\ldots,w_{K}\}$. The values of these parameters can
then be efficiently determined using a gradient-based algorithm; in
particular, we used Infinitesimal Perturbation Analysis (IPA) to evaluate the
objective function gradient as shown in \cite{cassandras2013optimal}.

Thus, a reasonable approach that has been suggested is to assign each agent a
linear trajectory. The 2D persistent monitoring problem would then be
formulated as consisting of the following tasks: $(i)$ finding $N$ linear
trajectories in terms of their length and exact location in $\Omega$, noting
that one or more agents may share one of these trajectories, and $(ii)$
controlling the motion of each agent on its trajectory. Task $(ii)$ is a
direct application of the 1D persistent monitoring problem solution, leaving
only task $(i)$ to be addressed. However, there is no reason to believe that a
linear trajectory is a good choice in a 2D setting. A broader choice is
provided by the set of elliptical trajectories which in fact encompass linear
ones when the minor axis of the ellipse becomes zero. Thus, we first proceed
with a comparison of these two types of trajectories. The main result of this
section is to formally show that an elliptical trajectory outperforms a linear
one using the average uncertainty metric in (\ref{eq:costfunction}) as the
basis for such comparison.

\begin{figure}[ptb]
\centering
\includegraphics[
height=1.34in,
width=3.27in]
{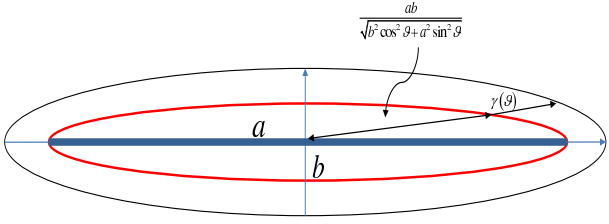}\caption{The red ellipse represents the agent trajectory. The
area defined by the black curve is the agent's effective coverage area.
$\frac{ab}{\sqrt{b^{2}cos^{2}(\vartheta)+a^{2}sin^{2}(\vartheta)}}%
+\gamma(\vartheta)$ is the distance between the ellipse center and the
coverage area boundary for a given $\vartheta$.}%
\label{effective area for proof}%
\end{figure}

To simplify notation, let $\omega=[x,y]\in\mathbb{R}^{2}$ and, for a single
agent, define
\begin{equation}
\Xi=\left\{  \omega\in\mathbb{R}^{2}|\exists t\in\lbrack0,T]\text{ such
that }Bp(\omega,s(t))>A(\omega)\right\}  \label{eq:effectiveArea}%
\end{equation}
Note that $\Xi$ above defines the \emph{effective coverage region} for the
agent, i.e., the region where the uncertainty corresponding to $R(\omega,t)$
with the dynamics in (\ref{eq:multiDynR}) can be strictly reduced given the
sensing capacity of the agent determined through $B$ and $p(\omega,s)$.
Clearly, $\Xi$ depends on the values of $s(t)$ which are dependent on the
agent trajectory. Let us define an elliptical trajectory so that the agent
position $s(t)=[s^{x}(t),s^{y}(t)]$ follows the general parametric form of an
ellipse:
\begin{equation}
\left\{
\begin{array}
[c]{c}%
s^{x}\left(  t\right)  =X+a\cos\rho\left(  t\right)  \cos\varphi-b\sin
\rho\left(  t\right)  \sin\varphi\\
s^{y}\left(  t\right)  =Y+a\cos\rho\left(  t\right)  \sin\varphi+b\sin
\rho\left(  t\right)  \cos\varphi
\end{array}
\right.  \label{ellipse}%
\end{equation}
where $[X,Y]$ is the center of the ellipse, $a,b$ are its major and minor axis
respectively, $\varphi\in\lbrack0,\pi)$ is the ellipse orientation (the angle
between the $x$ axis and the major ellipse axis) and $\rho(t)\in\lbrack
0,2\pi)$ is the eccentric anomaly of the ellipse. Assuming the agent moves
with constant maximal speed $1$ on this trajectory (based on (\ref{optu})), we
have $(\dot{s^{x}})^{2}+(\dot{s^{y}})^{2}=1$, which gives
\begin{align}
\dot{\rho}\left(  t\right) & =\left[(a\sin\rho(t)\cos\varphi+b\cos\rho(t)\sin
\varphi)^{2} \right.\nonumber \\
&\left. +(a\sin\rho(t)\sin\varphi-b\cos\rho(t)\cos\varphi)^{2}\right]^{-1/2}
\label{eq:speedEllipse}%
\end{align}
In order to make a fair comparison between a linear and an elliptical
trajectory, we normalize the objective function in (\ref{eq:costfunction})
with respect to the coverage area in (\ref{eq:effectiveArea}) and consider all
points in $\Xi$ (rather than discretizing it or limiting ourselves to a
finite set of sampling points). Thus, we define:
\begin{equation}
J(b)=\frac{1}{\Psi_{\Xi}}\int_{0}^{T}\int_{\Xi}R\left(  \omega,t\right)
d\omega dt \label{eq:fairCost}%
\end{equation}
where $\Psi_{\Xi}=\int_{\Xi}d\omega$ is the area of the effective
coverage region. Note that we view this normalized metric as a function of
$b\geq0$, so that when $b=0$ we obtain the uncertainty corresponding to a
linear trajectory. For simplicity, the trajectory is selected so that $[X,Y]$
coincides with the origin and $\varphi=0$, as illustrated in Fig.
\ref{effective area for proof} with the major axis $a$ assumed fixed.
Regarding the range of $b$, we will only be interested in values which are
limited to a neighborhood of zero that we will denote by $\mathcal{B}$. Given
$a$, this set dictates the values that $s(t)\in\Xi$ is allowed to take.
Finally, we make the following assumptions:

\textbf{Assumption 1}: $p(\omega,s)\equiv$ $p(D(\omega,s))$ is a continuous
function of $D(\omega,s)\equiv||\omega-s||$.

\textbf{Assumption 2}: Let $\omega,\omega^{\prime}$ be symmetric points in
$\Xi$ with respect to the center point of the ellipse. Then, $A(\omega
)=A(\omega^{\prime})$.

The first assumption simply requires that the sensing range of an agent is
continuous and the second that all points in $\Xi$ are treated uniformly
(as far as how uncertainty is measured) with respect to an elliptical
trajectory centered in this region.
The following result establishes the fact that an elliptical trajectory with
some $b>0$ can achieve a lower cost than a linear trajectory (i.e., $b=0$) in
terms of a long-term average uncertainty per unit area.

\begin{proposition}
\label{lem:ellipse} Under Assumptions 1-2 and $b\in\mathcal{B}$,
\[
\lim_{T\rightarrow\infty,b\rightarrow0}\frac{\partial J(b)}{\partial b}<0
\]
i.e., switching from a linear to an elliptical trajectory reduces the cost in
(\ref{eq:fairCost}).
\end{proposition}

\textbf{Proof}. Since a linear trajectory is the limit of an elliptical one
(with the major axis kept fixed) as the minor axis reaches $b=0$, our proof is
based on perturbing the minor axis $b$ away from $0$ and showing that we can
then achieve a lower average cost $J$ in (\ref{eq:fairCost}), as long as this
is measured over a sufficiently long time interval.

Obviously, the effective coverage area $\Psi_{\Xi}$ depends on the agent's
trajectory and, in particular, on the minor axis length $b$. From the
definition of $\Xi$ in (\ref{eq:effectiveArea}), note that $\Psi_{\Xi}$
monotonically increases in $b\in\mathcal{B}$, i.e., $\frac{\partial
\Psi_{\Xi}}{\partial b}>0$ and it immediately follows that:
\begin{equation}
\frac{\partial}{\partial b}(\frac{1}{\Psi_{\Xi}})=-\frac{\partial
\Psi_{\Xi}}{\partial b}\frac{1}{\Psi_{\Xi}^{2}}<0\label{eq:derivPSIwrtb}%
\end{equation}
We now rewrite the area integral in (\ref{eq:fairCost}) in a polar coordinate
system with $\omega=(\xi,\vartheta)\in\mathbb{R}^{2}$, where $\xi$ is the
polar radius and $\vartheta$ is the polar angle:
\begin{equation}
J(b)=\frac{1}{\Psi_{\Xi}}\int_{0}^{T}\int_{0}^{2\pi}\int_{0}%
^{G(a,b,\vartheta)+\gamma(\vartheta)}R(\xi,\vartheta,t)\xi d\xi d\vartheta
dt\label{eq:fairCost2}%
\end{equation}
where
\begin{equation}
G(a,b,\vartheta)=\frac{ab}{\sqrt{b^{2}\cos^{2}(\vartheta)+a^{2}\sin
^{2}(\vartheta)}}\label{G()}%
\end{equation}
is the ellipse equation in the polar coordinate system and $\gamma(\vartheta)$
is defined for any $(\xi,\vartheta)\in\mathbb{R}^{2}$ as
\begin{equation}
\gamma(\vartheta)=\sup_{\xi}\{Bp(\xi,\vartheta,s(t))>A(\xi,\vartheta
)\}-G(a,b,\vartheta)\label{eq:gammaTheta}%
\end{equation}
where $\sup_{\xi}\{Bp(\xi,\vartheta,s(t))>A(\xi,\vartheta)\}$ is the distance
between the ellipse center and the furthest point $(\xi,\vartheta)$, for any
given $\vartheta$, that can be effectively covered by the agent on the
ellipse.
Taking partial derivatives in (\ref{eq:fairCost2}) with respect to $b$, we
get
\begin{align}
\frac{\partial J}{\partial b} &  =-\frac{\partial\Psi_{\Xi}}{\partial
b}\frac{1}{\Psi_{\Xi}^{2}}\int_{0}^{T}\int_{\Xi}R\left(  \omega
,t\right)  d\omega dt\nonumber\\
&  +\frac{1}{\Psi_{\Xi}}\int_{0}^{T}\int_{0}^{2\pi}\left[  R\left(
G(a,b,\vartheta)+\gamma(\vartheta),\vartheta,t\right) \right.\nonumber \\
& \left. \cdot(G(a,b,\vartheta
)+\gamma(\vartheta))\cdot\frac{\partial G(a,b,\vartheta)}{\partial b}\right.
\nonumber\\
&  \left.  +\int_{0}^{G(a,b,\vartheta)+\gamma(\vartheta)}\frac{\partial
R\left(  \xi,\vartheta,t\right)  }{\partial b}\xi d\xi\right]  d\vartheta
dt\ \label{partialJwrtb}%
\end{align}
Recall that our objective is to show that when we perturb a linear trajectory
into an elliptical one, which is achieved by increasing $b$ from $0$ to some
small $b_{\epsilon}>0$, we can achieve a lower cost. Thus, we aim to show
$\frac{\partial J}{\partial b}|_{b\rightarrow0}<0$. From
(\ref{eq:derivPSIwrtb}), the first term of (\ref{partialJwrtb}) is negative,
therefore, we only need to show the second term is non-positive when
$b\rightarrow0$. By the definition (\ref{G()}), observe that when
$b\rightarrow0$, $G(a,b,\vartheta)\rightarrow0$, and $\frac{\partial G\left(
a,b,\vartheta\right)  }{\partial b}|_{b\rightarrow0}=\frac{1}{\sin\vartheta}$, for 
$\vartheta \neq 0$ and $\pi$; $\frac{\partial G\left(
a,b,\vartheta\right)  }{\partial b}|_{b\rightarrow0}=a$ for $\vartheta = 0$ or $\pi$.
Thus, the double integral of the second term of (\ref{partialJwrtb}) becomes
\begin{equation}
\int_{0}^{T}\int_{0}^{2\pi}\left[  \frac{\gamma(\vartheta)}{\sin\vartheta
}R\left(  \gamma(\vartheta),\vartheta,t\right)  +\int_{0}^{\gamma(\vartheta
)}\frac{\partial R\left(  \xi,\vartheta,t\right)  }{\partial b}\xi
d\xi\right]  d\vartheta dt\label{eq:partialJPartialBSecondTerm}%
\end{equation}
By Assumption 2, $A(\omega)=A(\omega^{\prime})$, where $\omega$ and
$\omega^{\prime}$ are symmetric with respect to the center point of the
ellipse, thus $A(\xi,\vartheta)=A(\xi,\vartheta+\pi)$. Then, for any
uncertainty value $R\left(  \gamma(\vartheta),\vartheta,t\right)  $ satisfying
(\ref{eq:multiDynR}), we can find $R\left(  \gamma(\vartheta+\pi
),\vartheta+\pi,t\right)  $ which is symmetric to it with respect to the
center point of the ellipse. Then, from (\ref{eq:gammaTheta}) and Fig.
\ref{effective area for proof}, note that $\gamma(\vartheta)=\gamma
(\vartheta+\pi)$. From the perspective of the point $(\gamma(\vartheta
),\vartheta)$, the agent movement observed with an initial position
$\rho(0)=\eta$ (following the dynamics in (\ref{eq:speedEllipse})) is the same
as the movement observed from $(\gamma(\vartheta+\pi),\vartheta+\pi)$ if the
agent starts from $\rho(0)=\eta+\pi$ when $T\rightarrow\infty$, since the cost
in (\ref{eq:fairCost}) is independent of initial conditions as $T\rightarrow
\infty$. Thus $R\left(  \gamma(\vartheta),\vartheta,t\right)  =R\left(
\gamma(\vartheta+\pi),\vartheta+\pi,t\right)  $. Since, in addition,
$\sin\vartheta=-\sin(\vartheta+\pi$), we have $\gamma(\vartheta)\frac
{R(\gamma(\vartheta),\vartheta,t)}{\sin\vartheta}=-\gamma(\vartheta+\pi
)\frac{R(\gamma(\vartheta+\pi),\vartheta+\pi,t)}{\sin(\vartheta+\pi)}$ and it
follows that
\begin{equation}
\lim_{T\rightarrow\infty,b\rightarrow0}\int_{0}^{T}\int_{0}^{2\pi}\frac
{\gamma(\vartheta)}{\sin\vartheta}R\left(  \gamma(\vartheta),\vartheta
,t\right)  d\vartheta dt=0\label{eq:partialJPartialBSecondTerm2}%
\end{equation}
We now turn our attention to the last integral of (\ref{partialJwrtb}). Two
cases need to be considered here in view of (\ref{eq:multiDynR}):

$(i)$ If $\exists t^{\prime}$ such that $R(\xi,\vartheta,t^{\prime})=0$ for
$t^{\prime}\in(0,t)$, then let
\begin{equation}
\tau_{f}(t)=\sup_{\tau\leq t}\{\tau:R(\xi,\vartheta,\tau)=0\} \label{tauf}%
\end{equation}
If $\tau_{f}(t)<t$, then $R(\xi,\vartheta,\tau)>0$ for all $\tau\in\lbrack
\tau_{f}(t),t)$ and $\tau_{f}(t)$ is the last time instant prior to $t$ when
$R\left(  \xi,\vartheta,\tau\right)  $ leaves an arc such that $R(\xi
,\vartheta,\tau)=0$. We can then write $R\left(  \xi,\vartheta,t\right)
=\int_{\tau_{f}(t)}^{t}\dot{R}\left(  \xi,\vartheta,\delta\right)  d\delta$.
Therefore,
\begin{align}
\frac{\partial R(\xi,\vartheta,t)}{\partial b}&=\frac{\partial t}{\partial
b}\dot{R}(\xi,\vartheta,t)-\frac{\partial\tau_{f}(t)}{\partial b}\dot{R}%
(\xi,\vartheta,\tau_{f}(t))\nonumber \\
&+\int_{\tau_{f}(t)}^{t}\frac{\partial\dot{R}%
(\xi,\vartheta,\delta)}{\partial b}d\delta\label{eq:partialRpartialB}%
\end{align}
Clearly, $\frac{\partial t}{\partial b}=0$ and since $\tau_{f}(t)$ is a time
instant when $R(\xi,\vartheta,t)$ leaves $R(\xi,\vartheta,t)=0$ then, by
Assumption 1, $\dot{R}(\xi,\vartheta,t)$ is a continuous function and we have
$\dot{R}(\xi,\vartheta,\tau_{f}(t))=0$. Therefore, (\ref{eq:partialRpartialB})
becomes
\begin{equation}
\frac{\partial R(\xi,\vartheta,t)}{\partial b}=\int_{\tau_{f}}^{t}%
\frac{\partial\dot{R}(\xi,\vartheta,\delta)}{\partial b}d\delta
\label{eq:partialJpartialY1}%
\end{equation}
where, from (\ref{eq:multiDynR}), $\dot{R}(\xi,\vartheta,\delta)=A(\xi
,\vartheta)-Bp(\xi,\vartheta,s(\delta))$.

If, on the other hand, $\tau_{f}(t)=t$, then $R(\xi,\vartheta,t)=0$ and we
define $\sigma_{f}(t)=\sup_{\sigma\leq t}\{\sigma:R(\xi,\vartheta,\sigma
)>0\}$. Proceeding as above, we get%
\[
\frac{\partial R(\xi,\vartheta,t)}{\partial b}=\int_{\sigma_{f}}^{t}%
\frac{\partial\dot{R}(\xi,\vartheta,\delta)}{\partial b}d\delta
\]
where now $\dot{R}(\xi,\vartheta,\delta)=0$ and we get%
\begin{equation}
\frac{\partial R(\xi,\vartheta,t)}{\partial b}=0 \label{partialRwrtb_zero}%
\end{equation}

$(ii)$ $R(\xi,\vartheta,t^{\prime})>0$ for all $t^{\prime}\in(0,t)$. In this
case, we define $\tau_{f}(t)=0$ and we have $R(\xi,\vartheta,t)=R(\xi
,\vartheta,0)+\int_{\tau_{f}(t)}^{t}\dot{R}\left(  \xi,\vartheta
,\delta\right)  d\delta$, where $\dot{R}\left(  \xi,\vartheta,\delta\right)
=A(\xi,\vartheta)-Bp(\xi,\vartheta,s(t)).$ Thus,
\begin{equation}
\frac{\partial R(\xi,\vartheta,t)}{\partial b} =\frac{\partial R(\xi
,\vartheta,0)}{\partial b}+\frac{\partial t}{\partial b}\dot{R}(\xi
,\vartheta,t) +\int_{\tau_{f}(t)}^{t}\frac{\partial\dot{R}(\xi,\vartheta
,\delta)}{\partial b}d\delta\label{eq:partialRpartialB0}%
\end{equation}
Clearly, $\frac{\partial t}{\partial b}=0$ and $\frac{\partial R(\xi
,\vartheta,0)}{\partial b}=0$, since $R(\xi,\vartheta,0)$ is the initial
uncertainty value at $(\xi,\vartheta)$ Then, (\ref{eq:partialRpartialB0})
becomes
\begin{equation}
\frac{\partial R(\xi,\vartheta,t)}{\partial b}=\int_{\tau_{f}}^{t}%
\frac{\partial\dot{R}(\xi,\vartheta,\delta)}{\partial b}d\delta
\label{eq:partialRpartialB01}%
\end{equation}
which is the same result as (\ref{eq:partialJpartialY1}).

Let us start by setting aside the much simpler case where
(\ref{partialRwrtb_zero}) applies and consider (\ref{eq:partialJpartialY1})
and (\ref{eq:partialRpartialB01}). Noting that $\frac{\partial A(\xi
,\vartheta)}{\partial b}=0$ we get
\begin{equation}
\frac{\partial\dot{R}\left(  \xi,\vartheta,\delta\right)  }{\partial
b}=-B\frac{\partial p(\xi,\vartheta,s(\delta))}{\partial b}
\label{eq:partialRpartialB2}%
\end{equation}
Recall that $[X,Y]$ has been selected to be the origin and that $\varphi=0$.
In this case, (\ref{ellipse}) becomes
\begin{equation}
s^{x}\left(  t\right)  =a\cos\rho\left(  t\right)  ,\text{ \ \ }s^{y}\left(
t\right)  =b\sin\rho\left(  t\right)  \label{eq:positionOnEllipse2}%
\end{equation}
Observing that $s^{x}(t)$ is independent of $b$, (\ref{eq:partialRpartialB2})
gives
\begin{align}
\frac{\partial\dot{R}\left(  \xi,\vartheta,\delta\right)  }{\partial b}  &
=-B\frac{\partial p(\xi,\vartheta,s(\delta))}{\partial s^{y}(\delta)}%
\frac{\partial s^{y}(\delta)}{\partial b}\nonumber\\
&  =-B\frac{\partial p(\xi,\vartheta,s(\delta))}{\partial D(\xi,\vartheta
,s(\delta))}\frac{\partial D(\xi,\vartheta,s(\delta))}{\partial s^{y}(\delta
)}\sin\rho(\delta) \label{eq:partialRpartialB3}%
\end{align}
where $D(\xi,\vartheta,s(\delta))=[(s^{x}(\delta)-\xi\cos\vartheta)^{2}%
+(s^{y}(\delta)-\xi\sin\vartheta)^{2}]^{1/2}$, hence
\begin{equation}
\frac{\partial D(\xi,\vartheta,s(\delta))}{\partial s^{y}(\delta)}=\frac
{s^{y}(\delta)-\xi\sin\vartheta}{D(\xi,\vartheta,s(\delta))}
\label{eq:partialRpartialB4}%
\end{equation}
Using (\ref{eq:partialRpartialB4}), (\ref{eq:partialRpartialB3}),
(\ref{eq:partialJpartialY1}) in the second integral of
(\ref{eq:partialJPartialBSecondTerm}), this integral becomes%
\begin{align}
&  \int_{0}^{T}\int_{0}^{2\pi}\int_{0}^{\gamma(\vartheta)}\frac{\partial
R\left(  \xi,\vartheta,t\right)  }{\partial b}\xi d\xi d\vartheta
dt\nonumber\\
&  =-B\int_{0}^{T}\int_{0}^{2\pi}\int_{0}^{\gamma(\vartheta)}\xi\int_{\tau
_{f}}^{t}\frac{\partial p(\xi,\vartheta,s(\delta))}{\partial D(\xi
,\vartheta,s(\delta))}\frac{\left(  s_{y}\left(  \delta\right)  -\xi
\sin\vartheta\right)  }{D(\xi,\vartheta,s(\delta))}\nonumber \\
&\cdot\sin\rho
(\delta)d\delta d\xi d\vartheta dt \label{eq:partialJParitialBSecondTerm3}%
\end{align}
Note that when $b\rightarrow0$, we have $s_{y}(\delta)\rightarrow0$. In
addition, $p(\xi,\vartheta,s(\delta))$ is a direct function of $D(\xi
,\vartheta,s(\delta))$, so that $\frac{\partial p(\xi,\vartheta,s(\delta
))}{\partial D(\xi,\vartheta,s(\delta))}$ is not an explicit function of
$\xi,\vartheta$ or $\delta$. Moreover, $\sin\rho(\delta)$ is not a function of
$\vartheta$. Thus, switching the integration order in
(\ref{eq:partialJParitialBSecondTerm3}) we get
\[
B\frac{\partial p(D)}{\partial D}\int_{0}^{T}\int_{\tau_{f}}^{t}\sin
\rho\left(  \delta\right)  \int_{0}^{2\pi}\int_{0}^{\gamma(\vartheta)}%
\frac{\xi^{2}\sin\vartheta}{D(\xi,\vartheta,s(\delta))}d\xi d\vartheta d\delta
dt
\]
Using Assumption 2, we make a symmetry argument similar to the one regarding
(\ref{eq:partialJPartialBSecondTerm2}). For any point $\omega=(\xi
,\vartheta)\in\mathbb{R}^{2}$, we can find $(\xi,\vartheta+\pi)$ which is
symmetric to it with respect to the center point of the ellipse and Assumption
2 implies that $A(\xi,\vartheta)=A(\xi,\vartheta+\pi)$. Then, from the
perspective of the point $(\xi,\vartheta)$, the agent movement observed with
an initial position $\rho(0)=\eta$ (following the dynamics in
(\ref{eq:speedEllipse})) is the same as the movement observed from
$(\xi,\vartheta+\pi)$ if the agent starts from $\rho(0)=\eta+\pi$ when
$T\rightarrow\infty$, since the cost in (\ref{eq:fairCost}) is independent of
initial conditions as $T\rightarrow\infty$. In addition, we again have
$\gamma(\vartheta)=\gamma(\vartheta+\pi)$, so that $\int_{0}^{\gamma
(\vartheta)}\frac{\sin\vartheta}{D(\xi,\vartheta,s(\delta))}=-\int_{0}%
^{\gamma(\vartheta+\pi)}\frac{\sin(\vartheta+\pi)}{D(\xi,\vartheta
+\pi,s(\delta))}$. Therefore,
\begin{equation}
\lim_{T\rightarrow\infty}\int_{0}^{2\pi}\int_{0}^{\gamma(\vartheta)}\frac
{\xi^{2}\sin\vartheta}{D(\xi,\vartheta,s(\delta))}d\xi d\vartheta=0
\end{equation}
and the second term of (\ref{eq:partialJPartialBSecondTerm}) gives
\begin{equation}
\lim_{T\rightarrow\infty,b\rightarrow0}\int_{0}^{T}\int_{0}^{2\pi}\int
_{0}^{\gamma(\vartheta)}\frac{\partial R\left(  \xi,\vartheta,t\right)
}{\partial b}\xi d\xi d\vartheta dt=0
\label{eq:partialJPartialBSecondTermZero}%
\end{equation}
In view of (\ref{eq:partialJPartialBSecondTerm2}) and
(\ref{eq:partialJPartialBSecondTermZero}), we have shown that the second term
of (\ref{partialJwrtb}) is $0$ and we are left with the first negative term
from (\ref{eq:derivPSIwrtb}), giving the desired result:
\begin{equation}
\lim_{T\rightarrow\infty,b\rightarrow0}\frac{\partial J(b)}{\partial b}%
=-\frac{\partial\Psi_{\Xi}}{\partial b}\frac{1}{\Psi_{\Xi}^{2}}\int
_{0}^{T}\int_{\Xi}R\left(  \omega,t\right)  d\omega dt<0
\label{eq:partialJPartialBNegative}%
\end{equation}
Finally, if (\ref{partialRwrtb_zero}) applies instead of
(\ref{eq:partialJpartialY1}), then (\ref{partialRwrtb_zero}) and
(\ref{eq:partialJPartialBSecondTerm2}) immediately imply that the second term
of (\ref{partialJwrtb}) is $0$, completing the proof. $\blacksquare$

Thus, we have proved that as $T\rightarrow\infty$, when $b$ is perturbed from
$0$ to some $b_{\epsilon}>0$, an elliptical trajectory achieves a lower cost
than a linear one. In other words, we have shown that elliptical trajectories
are more suitable for a 2D mission space in terms of achieving near-optimal
results in solving problem \textbf{P1}.

In other words, Prop. IV.1 shows that elliptical trajectories
are more suitable for a 2D mission space in terms of
achieving near-optimal results in solving problem \textbf{P1}.

\section{Optimal Elliptical Trajectories}

Based on our analysis thus far, we now tackle the problem of determining
optimal solutions within the class of elliptical trajectories. Our approach is
to associate with each agent an elliptical trajectory, parameterize each such
trajectory by its center, orientation and major and minor axes, and then solve
\textbf{P1} as a parametric optimization problem. Note that this includes the
possibility that two agents share the same trajectory if the solution to this
problem results in identical parameters for the associated ellipses. Choosing
elliptical trajectories, which are most likely suboptimal relative to a trajectory
obtained through a TPBVP solution of \textbf{P1}, offers several practical
advantages in addition to reduced computational complexity. Elliptical
trajectories induce a periodic structure to the agent movements which provides
predictability. As a result, it is also easier to handle issues related to
collision avoidance.

For an elliptical trajectory, the $n$th agent movement is described as in
(\ref{ellipse}) by
\begin{equation}
\left\{
\begin{array}
[c]{c}%
s_{n}^{x}\left(  t\right)  =X_{n}+a_{n}\cos\rho_{n}\left(  t\right)
\cos\varphi_{n}-b_{n}\sin\rho_{n}\left(  t\right)  \sin\varphi_{n}\\
s_{n}^{y}\left(  t\right)  =Y_{n}+a_{n}\cos\rho_{n}\left(  t\right)
\sin\varphi_{n}+b_{n}\sin\rho_{n}\left(  t\right)  \cos\varphi_{n}%
\end{array}
\right.  \label{eq:positionOnEllipse}%
\end{equation}
where $[X_{n},Y_{n}]$ is the center of the $n$th ellipse, $a_{n},b_{n}$ are
its major and minor axes respectively and $\varphi_{n}\in\lbrack0,\pi)$ is its
orientation, i.e., the angle between the horizontal axis and the major axis of the
$n$th ellipse. Note that the parameter $\rho_{n}(t)\in\lbrack0,2\pi)$ is the
eccentric anomaly. Therefore, we replace problem \textbf{P1} by the
determination of optimal parameter vectors $\Upsilon_{n}\equiv\lbrack
X_{n},Y_{n},a_{n},b_{n},\varphi_{n}]^{\mathtt{T}},n=1,\dots,N$, and formulate
the following problem \textbf{P2}:%
\begin{equation}
\mathbf{P2}:\text{ \ \ }\min_{\Upsilon_{n},n=1,\ldots,N}\text{ \ }J=\int
_{0}^{T}\sum_{i=1}^{M}R_{i}(\Upsilon_{1},\ldots,\Upsilon
_{N},t)dt\label{P2}%
\end{equation}
Observe that the behavior of each agent under the optimal ellipse control
policy is that of a \emph{hybrid system} whose dynamics undergo switches when
$R_{i}(t)$ reaches or leaves the boundary value $R_{i}=0$ (the
\textquotedblleft events\textquotedblright\ causing the switches). As a
result, we are faced with a parametric optimization problem for a system with
hybrid dynamics. We solve this hybrid system problem using a gradient-based
approach in which we apply IPA to determine the gradients $\nabla
R_{i}(\Upsilon_{1},\ldots,\Upsilon_{N},t)$ on line (hence, $\nabla J$), i.e.,
directly using information from the agent trajectories and iterate upon them.

\subsection{Infinitesimal Perturbation Analysis (IPA)}

We begin with a brief review of the IPA framework for general stochastic
hybrid systems as presented in \cite{cassandras2009perturbation}. The purpose
of IPA is to study the behavior of a hybrid system state as a function of a
parameter vector $\theta\in\Theta$ for a given compact, convex set
$\Theta\subset\mathbb{R}^{l}$. Let $\{\tau_{k}(\theta)\}$, $k=1,\ldots,K$,
denote the occurrence times of all events in the state trajectory. For
convenience, we set $\tau_{0}=0$ and $\tau_{K+1}=T$. Over an interval
$[\tau_{k}(\theta),\tau_{k+1}(\theta))$, the system is at some mode during
which the time-driven state satisfies $\dot{x}\ =\ f_{k}(x,\theta,t)$. An
event at $\tau_{k}$ is classified as $(i)$ \emph{Exogenous} if it causes a
discrete state transition independent of $\theta$ and satisfies $\frac
{d\tau_{k}}{d\theta}=0$; $(ii)$ \emph{Endogenous}, if there exists a
continuously differentiable function $g_{k}:\mathbb{R}^{n}\times
\Theta\rightarrow\mathbb{R}$ such that $\tau_{k}\ =\ \min\{t>\tau
_{k-1}\ :\ g_{k}\left(  x\left(  \theta,t\right)  ,\theta\right)  =0\}$; and
$(iii)$ \emph{Induced} if it is triggered by the occurrence of another event
at time $\tau_{m}\leq\tau_{k}$. IPA specifies how changes in $\theta$
influence the state $x(\theta,t)$ and the event times $\tau_{k}(\theta)$ and,
ultimately, how they influence interesting performance metrics which are
generally expressed in terms of these variables.

We define:
\[
x^{\prime}(t)\equiv\frac{\partial x(\theta,t)}{\partial\theta},\text{ }%
\tau_{k}^{\prime}\equiv\frac{\partial\tau_{k}(\theta)}{\partial\theta},\text{
\ \ }k=1,\ldots,K
\]
for all state and event time derivatives. It is shown in
\cite{cassandras2009perturbation} that $x^{\prime}(t)$ satisfies:%
\begin{equation}
\frac{d}{dt}x^{\prime}\left(  t\right)  =\frac{\partial f_{k}\left(  t\right)
}{\partial x}x^{\prime}\left(  t\right)  +\frac{\partial f_{k}\left(
t\right)  }{\partial\theta} \label{eq:xpDyn}%
\end{equation}
for $t\in\lbrack\tau_{k},\tau_{k+1})$ with boundary condition:
\begin{equation}
x^{\prime}(\tau_{k}^{+})\ =\ x^{\prime}(\tau_{k}^{-})+\left[  f_{k-1}(\tau
_{k}^{-})-f_{k}(\tau_{k}^{+})\right]  \tau_{k}^{\prime} \label{eq:xpBoundary}%
\end{equation}
for $k=0,\ldots,K$,
where $\tau_{k}^{-}$ is the left limit of $\tau_{k}$.
In addition, in (\ref{eq:xpBoundary}), the gradient vector
for each $\tau_{k}$ is $\tau_{k}^{\prime}=0$ if the event at $\tau_{k}$ is
exogenous and
\begin{equation}
\tau_{k}^{\prime}=-\left[  \frac{\partial g_{k}}{\partial x}f_{k}(\tau_{k}%
^{-})\right]  ^{-1}\left(  \frac{\partial g_{k}}{\partial\theta}%
+\frac{\partial g_{k}}{\partial x}x^{\prime}(\tau_{k}^{-})\right)
\label{eq:taukp}%
\end{equation}
if the event at $\tau_{k}$ is endogenous (i.e., $\ g_{k}\left(  x\left(
\theta,\tau_{k}\right)  ,\theta\right)  =0$) and defined as long as
$\frac{\partial g_{k}}{\partial x}f_{k}(\tau_{k}^{-})\neq0$.

In our case, the parameter vectors are $\Upsilon_{n}\equiv\lbrack X_{n}%
,Y_{n},a_{n},b_{n},\varphi_{n}]^{\mathtt{T}}$ as defined earlier, and we seek
to determine optimal vectors $\varUpsilon_{n}^{\ast},$ $n=1,\dots,N$. We will
use IPA to evaluate $\nabla J(\varUpsilon_{1},\ldots,\varUpsilon_{N}%
)=[\frac{\partial J}{\partial \varUpsilon_{1}},\ldots,\frac{\partial J}{\partial \varUpsilon_{N}}%
]^{\mathtt{T}}$. From (\ref{P2}), this gradient clearly depends on $\nabla
R_{i}(t)=\left[  \frac{\partial R_{i}\left(  t\right)  }{\partial
\varUpsilon_{1}},\ldots,\frac{\partial R_{i}\left(  t\right)  }{\partial
\varUpsilon_{N}}\right]  ^{\mathtt{T}}$. In turn, this gradient depends on
whether the dynamics of $R_{i}(t)$ in (\ref{eq:multiDynR}) are given by
$\dot{R}_{i}(t)=0$ or $\dot{R}_{i}(t)=A_{i}-BP_{i}\left(  \mathbf{s}%
(t)\right)  $. The dynamics switch at event times $\tau_{k}$, $k=1,\ldots,K,$
when $R_{i}(t)$ reaches or escapes from $0$ which are observed on a trajectory
over $[0,T]$ based on a given $\varUpsilon_{n},$ $n=1,\dots,N$.

\textbf{IPA equations}. We begin by recalling the dynamics of $R_{i}\left(
t\right)  $ in (\ref{eq:multiDynR}) which depend on the relative positions of
all agents with respect to $[\alpha_{i},\beta_{i}]$ and change at time
instants $\tau_{k}$ such that either $R_{i}(\tau_{k})=0$ with $R_{i}(\tau
_{k}^{-})>0$ or $A_{i}>BP_{i}\left(  \mathbf{s}(\tau_{k})\right)  $ with
$R_{i}(\tau_{k}^{-})=0$. Moreover, the agent positions $s_{n}\left(  t\right)
=[s_{n}^{x}(t),s_{n}^{y}(t)]$, $n=1,\ldots,N$, on an elliptical trajectory are
expressed using (\ref{eq:positionOnEllipse}). Viewed as a hybrid system, we
can now concentrate on all events causing transitions in the dynamics of
$R_{i}\left(  t\right)  $, $i=1,\ldots,M$, since any other event has no effect
on the values of $\nabla R_{i}(\Upsilon_{1},\ldots,\Upsilon_{N},t)$ at
$t=\tau_{k}$.

For notational simplicity, we define $\omega_{i}=[\alpha_{i},\beta_{i}%
]\in\Omega$. First, if $R_{i}(t)=0$ and $A(\omega_{i})-BP(\omega
_{i},\mathbf{s}(t))\leq0$, applying (\ref{eq:xpDyn}) to $R_{i}\left(  t\right)  $
and using (\ref{eq:multiDynR}) gives
\begin{equation}
\frac{d}{dt}\frac{\partial R_{i}(t)}{\partial\varUpsilon_{n}}=0
\label{eq:pRpGamma1}%
\end{equation}
When $R_{i}(t)>0$, we have
\begin{equation}
\frac{d}{dt}\frac{\partial R_{i}\left(  t\right)  }{\partial\Upsilon_{n}%
}=-B\frac{\partial p_{n}\left(  \omega_{i},s_{n}\left(  t\right)  \right)
}{\partial\Upsilon_{n}}\textstyle\prod\limits_{d\neq n}^{N}\left[
1-p_{d}\left(  \omega_{i},s_{d}\left(  t\right)  \right)  \right]
\label{eq:pRpGamma}%
\end{equation}
Noting that
$p_{n}(\omega_{i},s_{n}(t))\equiv p_{n}(D(\omega_{i},s_{n}(t)))$, we have
\begin{equation}
\frac{\partial p_{n}\left(  \omega_{i},s_{n}\left(  t\right)  \right)
}{\partial\varUpsilon_{n}}=\frac{\partial p_{n}(D(\omega_{i},s_{n}(t)))}{\partial D(\omega_{i},s_{n}(t)))}\frac{\partial D(\omega_{i}%
,s_{n}(t))}{\partial\varUpsilon_{n}} \label{eq:pPpGamma}%
\end{equation}
where $D(\omega_{i},s_{n}(t))=[(s_{n}^{x}(t)-\alpha_{i})^{2}+(s_{n}%
^{y}(t)-\beta_{i})^{2}]^{1/2}$. For simplicity, we write $D=D(\omega_{i}%
,s_{n}(t))$ and we get
\begin{equation}
\frac{\partial D}{\partial\varUpsilon_{n}}=\frac{1}{2D}\left(  \frac{\partial
D}{\partial s_{n}^{x}}\frac{\partial s_{n}^{x}}{\partial\varUpsilon_{n}}%
+\frac{\partial D}{\partial s_{n}^{y}}\frac{\partial s_{n}^{y}}{\partial
\varUpsilon_{n}}\right)  \label{eq:pDpGamma}%
\end{equation}
where $\frac{\partial D}{\partial s_{n}^{x}}=2(s_{n}^{x}-\alpha_{i})$ and
$\frac{\partial D}{\partial s_{n}^{y}}=2(s_{n}^{y}-\beta_{i})$. Note that
$\frac{\partial s_{n}^{x}}{\partial\varUpsilon_{n}}=[\frac{\partial s_{n}^{x}%
}{\partial X_{n}},\frac{\partial s_{n}^{x}}{\partial Y_{n}},\frac{\partial
s_{n}^{x}}{\partial a_{n}},\frac{\partial s_{n}^{x}}{\partial b_{n}}%
,\frac{\partial s_{n}^{x}}{\partial\varphi_{n}}]^{\mathtt{T}}$ and
$\frac{\partial s_{n}^{y}}{\partial\varUpsilon_{n}}=[\frac{\partial s_{n}^{y}%
}{\partial X_{n}},\frac{\partial s_{n}^{y}}{\partial Y_{n}},\frac{\partial
s_{n}^{y}}{\partial a_{n}},\frac{\partial s_{n}^{y}}{\partial b_{n}}%
,\frac{\partial s_{n}^{y}}{\partial\varphi_{n}}]^{\mathtt{T}}$. From
(\ref{eq:positionOnEllipse}), for $\frac{\partial s_{n}^{x}}{\partial
\varUpsilon_{n}}$, we obtain%
\begin{align*}
\frac{\partial s_{n}^{x}}{\partial X_{n}}  &  =1,\text{ \ \ }\frac{\partial
s_{n}^{x}}{\partial Y_{n}}=0\\
\frac{\partial s_{n}^{x}}{\partial a_{n}}  &  =\cos\rho_{n}\left(  t\right)
\cos\varphi_{n},\text{ \ \ }\frac{\partial s_{n}^{x}}{\partial b_{n}}%
=-\sin\rho_{n}\left(  t\right)  \sin\varphi_{n}\\
\frac{\partial s_{n}^{x}}{\partial\varphi_{n}}  &  =-a_{n}\cos\rho_{n}\left(
t\right)  \sin\varphi_{n}-b\sin\rho_{n}\left(  t\right)  \cos\varphi_{n}%
\end{align*}
Similarly, for $\frac{\partial s_{n}^{y}}{\partial\varUpsilon_{n}}$, we get
$\frac{\partial s_{n}^{y}}{\partial X_{n}}=0,\frac{\partial s_{n}^{y}%
}{\partial Y_{n}}=1,\frac{\partial s_{n}^{y}}{\partial a_{n}}=\cos\rho
_{n}\left(  t\right)  \sin\varphi_{n},\frac{\partial s_{n}^{y}}{\partial
b_{n}}=\sin\rho_{n}\left(  t\right)  \cos\varphi_{n}$ and $\frac{\partial
s_{n}^{y}}{\partial\varphi_{n}}=a_{n}\cos\rho_{n}\left(  t\right)  \cos
\varphi_{n}-b\sin\rho_{n}\left(  t\right)  \sin\varphi_{n}.$ Using
$\frac{\partial s_{n}^{x}}{\partial\varUpsilon_{n}}$ and $\frac{\partial
s_{n}^{y}}{\partial\varUpsilon_{n}}$ in (\ref{eq:pDpGamma}) and then
(\ref{eq:pPpGamma}) and back into (\ref{eq:pRpGamma}), we can finally obtain
$\frac{\partial R_{i}\left(  t\right)  }{\partial\Upsilon_{n}}$ for
$t\in\lbrack\tau_{k},\tau_{k+1})$ as
\begin{equation}
\frac{\partial R_{i}\left(  t\right)  }{\partial\Upsilon_{n}}=\frac{\partial
R_{i}\left(  \tau_{k}^{+}\right)  }{\partial\Upsilon_{n}}+\left\{
\begin{array}
[c]{cc}%
0 & \text{if }R_{i}\left(  t\right)  =0,\\ & \text{ }A_{i}\leq BP_{i}\left(
\mathbf{s}\left(  t\right)  \right) \\
\int_{\tau_{k}}^{t}\frac{d}{dt}\frac{\partial R_{i}\left(  t\right)
}{\partial\Upsilon_{n}}dt & \text{otherwise}%
\end{array}
\right.  \label{eq:pRt}%
\end{equation}
where the integral above is obtained from (\ref{eq:pRpGamma1}%
)-(\ref{eq:pPpGamma}). Thus, it remains to determine the components $\nabla
R_{i}(\tau_{k}^{+})$ in (\ref{eq:pRt}) using (\ref{eq:xpBoundary}). This
involves the event time gradient vectors $\nabla\tau_{k}=\frac{\partial
\tau_{k}}{\partial\varUpsilon_{n}}$ for $k=1,\ldots,K$, which will be
determined through (\ref{eq:taukp}). There are two possible cases regarding
the events that cause switches in the dynamics of $R_{i}\left(  t\right)  $:

\emph{Case 1}: At $\tau_{k}$, $\dot{R}_{i}\left(  t\right)  $ switches from
$\dot{R}_{i}\left(  t\right)  =0$ to $\dot{R}_{i}\left(  t\right)
=A_{i}-BP_{i}\left(  \mathbf{s}(t)\right)  $. In this case, it is easy to see
that the dynamics $R_{i}(t)$ are continuous, so that $f_{k-1}(\tau_{k}%
^{-})=f_{k}(\tau_{k}^{+})$ in (\ref{eq:xpBoundary}) applied to $R_{i}(t)\ $and
we get%
\begin{equation}
\nabla R_{i}(\tau_{k}^{+})=\nabla R_{i}(\tau_{k}^{-}),\text{ }i=1,\ldots,M
\label{eq:wRi}%
\end{equation}

\emph{Case 2}: At $\tau_{k}$, $\dot{R}_{i}\left(  t\right)  $ switches from
$\dot{R}_{i}\left(  t\right)  =A_{i}-BP_{i}\left(  \mathbf{s}(t)\right)  $ to
$\dot{R}_{i}\left(  t\right)  =0$, i.e., $R_{i}(\tau_{k})$ becomes zero. In
this case, we need to first evaluate $\nabla\tau_{k}$ from (\ref{eq:taukp}) in
order to determine $\nabla R_{i}(\tau_{k}^{+})$ through (\ref{eq:xpBoundary}).
Observing that this event is endogenous, (\ref{eq:taukp}) applies with
$g_{k}=R_{i}=0$ and we get
\begin{equation}
\nabla\tau_{k}=-\frac{\nabla R_{i}(\tau_{k}^{-})}{A(\omega_{i})-BP(\omega
_{i},\mathbf{s}(\tau_{k}^{-}))}%
\end{equation}
It follows from (\ref{eq:xpBoundary}) that%
\begin{equation}
\nabla R_{i}(\tau_{k}^{+})=\nabla R_{i}(\tau_{k}^{-})-\frac{[A(\omega
_{i})-BP(\omega_{i},\mathbf{s}(\tau_{k}^{-}))]\nabla R_{i}\left(  \tau_{k}%
^{-}\right)  }{A(\omega_{i})-BP(\omega_{i},\mathbf{s}(\tau_{k}^{-}))}=0
\label{eq:wRpR0}%
\end{equation}
Thus, $\nabla R_{i}(\tau_{k}^{+})$ is always reset to $0$ regardless of
$\nabla R_{i}(\tau_{k}^{-})$.

\textbf{Objective Function Gradient Evaluation.} Based on our analysis, we
first rewrite $J$ in (\ref{P2}) as%
\[
J(\varUpsilon_{1},\ldots,\varUpsilon_{N})=\sum_{i=1}^{M}\sum_{k=0}%
^{K}{\displaystyle\int_{\tau_{k}(\varUpsilon_{1},\ldots,\varUpsilon_{N}%
)}^{\tau_{k+1}(\varUpsilon_{1},\ldots,\varUpsilon_{N})}}R_{i}(\varUpsilon_{1}%
,\ldots,\varUpsilon_{N},t)dt
\]
and (omitting some function arguments) we get%
\[
\nabla J={\displaystyle\sum_{i=1}^{M}}{\displaystyle\sum_{k=0}^{K}}\left(
{\displaystyle\int_{\tau_{k}}^{\tau_{k+1}}}\nabla R_{i}\left(  t\right)
dt+R_{i}\left(  \tau_{k+1}\right)  \nabla\tau_{k+1}-R_{i}\left(  \tau
_{k}\right)  \nabla\tau_{k}\right)
\]
Observing the cancelation of all terms of the form $R_{i}\left(  \tau
_{k}\right)  \nabla\tau_{k}$ for all $k$ (with $\tau_{0}=0$, $\tau_{K+1}=T$
fixed), we finally get
\begin{equation}
\nabla J(\varUpsilon_{1},\ldots,\varUpsilon_{N})=\sum_{i=1}^{M}%
\sum_{k=0}^{K}{\displaystyle\int_{\tau_{k}}^{\tau_{k+1}}}\nabla R_{i}\left(
t\right)  dt \label{eq:derivTotal}%
\end{equation}
This depends entirely on $\nabla R_{i}\left(  t\right)  $, which is obtained
from (\ref{eq:pRt}) and the event times $\tau_{k}$, $k=1,\ldots,K$, given
initial conditions $s_{n}\left(  0\right)  $ for $n=1,\ldots,N$, and
$R_{i}\left(  0\right)  $ for $i=1,\ldots,M$. In (\ref{eq:pRt}),
$\frac{\partial R_{i}\left(  \tau_{k}^{+}\right)  }{\partial\Upsilon_{n}}$ is
obtained through (\ref{eq:wRi})-(\ref{eq:wRpR0}), whereas $\frac{d}{dt}%
\frac{\partial R_{i}\left(  t\right)  }{\partial\Upsilon_{n}}$ is obtained
through (\ref{eq:pRpGamma1})-(\ref{eq:pDpGamma}).

\textbf{Remark 2}. Observe that the evaluation of $\nabla R_{i}\left(
t\right)  $, hence $\nabla J$, is \emph{independent} of $A_{i}$,
$i=1,\ldots,M$, i.e., the values in our uncertainty model. In fact, the
dependence of $\nabla R_{i}\left(  t\right)  $ on $A_{i}$, $i=1,\ldots,M$,
manifests itself through the event times $\tau_{k}$, $k=1,\ldots,K$, that do
affect this evaluation, but they, unlike $A_{i}$ which may be unknown, are
directly observable during the gradient evaluation process. Thus, the IPA
approach possesses an inherent \emph{robustness} property: there is no need to
explicitly model how uncertainty affects $R_{i}(t)$ in (\ref{eq:multiDynR}).
Consequently, we may treat $A_{i}$ as unknown without affecting the solution
approach (the values of $\nabla R_{i}\left(  t\right)  $ are obviously
affected). We may also allow this uncertainty to be modeled through random
processes $\{A_{i}(t)\}$, $i=1,\ldots,M$; in this case, however, the result of
Proposition \ref{lem:ellipse} no longer applies without some conditions on the
statistical characteristics of $\{A_{i}(t)\}$ and the resulting $\nabla J$ is
an \emph{estimate} of a stochastic gradient.

\textbf{Remark 3}. Note that the number of agents affects the number of
derivative components in (\ref{eq:derivTotal}), so the complexity of $\nabla
J(\varUpsilon_{1},\ldots,\varUpsilon_{N})$ in (\ref{eq:derivTotal}) grows
linearly in the number of agents $N$. In addition, the calculation of $\nabla
J(\varUpsilon_{1},\ldots,\varUpsilon_{N})$ in (\ref{eq:derivTotal}) grows
linearly in $T$, as a longer operation time only implies more events at whose
occurrence times $\tau_{k}$ the objective function gradient is updated. In
other words, solving the problem using IPA is scalable with respect to the
number of agents and the operation time.

\subsection{Objective Function Optimization}

We now seek to obtain $[\varUpsilon_{1}^{\ast},\ldots,\varUpsilon_{N}^{\ast}]$
minimizing $J(\varUpsilon_{1},\ldots,\varUpsilon_{N})$ through a standard
gradient-based optimization algorithm of the form%
\begin{equation}
\lbrack\varUpsilon_{1}^{l+1},\ldots,\varUpsilon_{N}^{l+1}]=[\varUpsilon_{1}%
^{l},\ldots,\varUpsilon_{N}^{l}]-[\eta_{1}^{l},\ldots,\eta_{N}^{l}%
]\tilde{\nabla}J(\varUpsilon_{1}^{l},\ldots,\varUpsilon_{N}^{l})
\label{eq:updatetheta}%
\end{equation}
where $\{\eta_{n}^{l}\}$, $l=1,2,\ldots$ are appropriate step size sequences
and $\tilde{\nabla}J(\varUpsilon_{1}^{l},\ldots,\varUpsilon_{N}^{l})$ is the
projection of the gradient $\nabla J(\varUpsilon_{1},\ldots,\varUpsilon_{N})$
onto the feasible set, i.e., $s_{n}(t)\in\Omega$ for all $t\in\lbrack0,T]$,
$n=1,\ldots,N$.
The optimization algorithm terminates when $|\tilde{\nabla
}J(\varUpsilon_{1}^{l},\ldots,\varUpsilon_{N}^{l})|<\varepsilon$ (for a fixed
threshold $\varepsilon$) for some $[\varUpsilon_{1}^{\ast},\ldots
,\varUpsilon_{N}^{\ast}]$. When $\varepsilon>0$ is small, $[\varUpsilon_{1}^{l},\ldots,\varUpsilon_{N}^{l}]$ is
believed to be in the neighborhood of the local optimum, then we set
$[\varUpsilon_{1}^{\ast},\ldots,\varUpsilon_{N}^{\ast}]=[\varUpsilon_{1}^{l}%
,\ldots,\varUpsilon_{N}^{l}].$ However, in our problem the function
$J(\varUpsilon_{1},\ldots,\varUpsilon_{N})$ is non-convex and there are
actually many local optima depending on the initial controllable parameter
vector $[\varUpsilon_{1}^{0},\ldots,\varUpsilon_{N}^{0}]$. In the next
section, we propose a stochastic comparison algorithm which addresses this
issue by randomizing over the initial points $[\varUpsilon_{1}^{0}%
,\ldots,\varUpsilon_{N}^{0}]$ . This algorithm defines a process which
converges to a global optimum under certain well-defined conditions.

\section{Stochastic Comparison Algorithm for global optimality}

Gradient-based optimization algorithms are generally efficient and effective
in finding the global optimum when one is uniquely specified by the point
where the gradient is zero. When this is not the case, to seek a global
optimum one must resort to several alternatives which include a variety of
random search algorithms. In this section, we use the Stochastic Comparison
algorithm in \cite{bao1996stochastic} to find the global optimum. As shown in
\cite{bao1996stochastic}, for a stochastic system, if $(i),$ the cost function
$J(\varUpsilon)$ is continuous in $\varUpsilon$ and $(ii),$ for each estimate
$\hat{J}(\varUpsilon)$ of $J(\varUpsilon)$ the error $W(\varUpsilon)=\hat
{J}(\varUpsilon)-J(\varUpsilon)$ has a symmetric pdf, then the Markov process
$\{\varUpsilon_{k}\}$ generated by the Stochastic Comparison algorithm will
converge to an $\epsilon-$optimal interval of the global optimum for
arbitrarily small $\epsilon>0.$ In short,
$\lim_{k\rightarrow\infty}P[\varUpsilon^{k}\in\varUpsilon_{\epsilon}^{\ast
}]=1,\text{for any }\epsilon>0$,
where $\varUpsilon_{\epsilon}^{\ast}$ is defined as $\varUpsilon_{\epsilon
}^{\ast}=\{\varUpsilon|J(\varUpsilon)\leq J(\varUpsilon^{\ast})+\epsilon\}$.
Using the Continuous Stochastic Comparison (CSC) Algorithm developed in
\cite{bao1996stochastic} for a general continuous optimization problem,
consider $\varUpsilon\in\Phi$ to be a controllable vector, where $\Phi$ is the
bounded feasible controllable parameter space. The Stochastic Comparison
Algorithm is presented in \textbf{Algorithm 1}. \begin{algorithm}
\caption{: Continuous Stochastic Comparison (CSC) Algorithm.}
\label{alg:DSC}
\begin{algorithmic}[1]
\STATE Initialize $\varUpsilon^{0} = \phi^{0},k=0.$
\STATE For a given $\varUpsilon^{k} = \phi^{k}$, sample the next candidate point $Z^{k}$ from $\Phi$ according to a uniform distribution over $\Phi.$
\STATE For a given $Z^{k} = \zeta^{k}$, set
\begin{equation}
\Upsilon ^{k+1}=\left\{
\begin{array}{cc}
Z^{k}, & \text{with probability }p^{k}\text{,} \\
\Upsilon ^{k}, & \text{with probability }1-p^{k}\text{,}%
\end{array}%
\right.
\end{equation}
where $p^{k} = \{P[\hat{J}(\zeta^{k})<\hat{J}(\phi^{k})]\}^{L_{k}}.$
\STATE Replace $k$ by $k+1,$ and go to Step 2.
\end{algorithmic}
\end{algorithm}In the CSC algorithm, the probability $p^{k}$ is actually not
calculable, since we do not know the underlying probability functions.
However, it is realizable in the following way: both $\hat{J}(\zeta^{k})$ and
$\hat{J}(\phi^{k})$ are estimated $L_{k}$ times for an appropriately selected
increasing sequence $\{L_{k}\}$. If $\hat{J}(\zeta^{k})<\hat{J}(\phi^{k})$
every time, we set $\varUpsilon^{k+1}=Z^{k}.$ Otherwise, we set
$\varUpsilon^{k+1}=\varUpsilon^{k}.$

As discussed in \textbf{Remark 3}, the persistent monitoring problem
\textbf{P2} becomes a stochastic optimization problem if $A_{i}(t)$,
$i=1,\ldots,M$, are stochastic processes. However, for the deterministic
setting in which all $A_{i}$ are constant, the observed value $\hat{J}$
coincides with the actual value $J$ and a one-time comparison $\hat{J}%
(\zeta^{k})<\hat{J}(\phi^{k})$ is sufficient to replace $\phi^{k}$ with
$\zeta^{k}$ for $\varUpsilon^{k+1}.$ In this case, step 3 in \textbf{Algorithm
1} becomes, for a given $Z^{k}=\zeta^{k}$:
\begin{equation}
\varUpsilon^{k+1}=\left\{
\begin{array}
[c]{cc}%
Z^{k} & \text{if }J(\zeta^{k})<J(\phi^{k})\\
\varUpsilon^{k} & \text{otherwise}%
\end{array}
\right.
\end{equation}
and the CSC algorithm in this deterministic setting reduces to a comparison
algorithm with multi-starts over the 6-dimensional controllable vector
$\Upsilon_{n}\equiv\lbrack X_{n},Y_{n},a_{n},b_{n},\varphi_{n},\rho
_{n}]^{\mathtt{T}},$ for each ellipse associated with agent $n=1,\dots,N$.

\begin{algorithm}
\caption{: IPA-based Optimization Algorithm using CSC to find $\Upsilon_{n}$, $n=1,\ldots,N.$}
\label{alg:IPA}
\begin{algorithmic}[1]
\STATE Set $\epsilon>0$, $k=0$. Initialize $\varUpsilon^{0} = \phi^{0}$, 
where $\phi^{0} = [\varUpsilon_{1}^{0},\ldots,\varUpsilon_{N}^{0}].$
Initialize $L_{0}$, where $\{L_{k}\}$ is an appropriately selected increasing sequence.
\WHILE {$k < K,$}
\STATE For a given $\varUpsilon^{k} = \phi^{k}$,
\REPEAT
\STATE Compute $s_{n}(t)$,  $t\in[0,T]$ using (\ref{eq:positionOnEllipse}) and $\phi^{k}$ for $n=1,\ldots,N$
\STATE Compute $\hat J(\phi^{k})$, $\tilde\nabla J(\phi^{k})$ and update $\phi^{k}$ through \eqref{eq:updatetheta}.
\UNTIL{$|\tilde\nabla J(\phi^{k})|<\epsilon$}
\STATE Sample the next candidate point $Z^{k}$ from $\Phi$ according to a uniform distribution over $\Phi.$ For a given $Z^{k} = \zeta^{k}$,
\REPEAT
\STATE Compute $s_{n}(t)$,  $t\in[0,T]$ using (\ref{eq:positionOnEllipse}) and $\zeta^{k}$ for $n=1,\ldots,N$
\STATE Compute $\hat J(\zeta^{k})$, $\tilde\nabla J(\zeta^{k})$ and update $\zeta^{k}$ through \eqref{eq:updatetheta}.
\UNTIL{$|\tilde\nabla J(\zeta^{k})|<\epsilon$}
\STATE Set
\begin{equation}
\Upsilon ^{k+1}=\left\{
\begin{array}{cc}
Z^{k}, & \text{with probability }p^{k}\text{,} \\
\Upsilon ^{k}, & \text{with probability }1-p^{k}\text{,}%
\end{array}%
\right.
\end{equation}
where $p^{k} = \{P[\hat{J}(\zeta^{k})<\hat{J}(\phi^{k})]\}^{L_{k}}.$
\STATE Replace $k$ by $k+1$.
\ENDWHILE
\STATE Set $\varUpsilon^{\ast} = \varUpsilon^{K}$.

\end{algorithmic}
\end{algorithm}

\section{Numerical Results}

\label{sec:numericalResults} 
\begin{figure}[ptb]
\centering
\subfigure[Red ellipses are the initial trajectories and blue ellipses are the final trajectories.]{
\label{fig:T200Tra}
\includegraphics[
height=1.8in,
width=3in]
{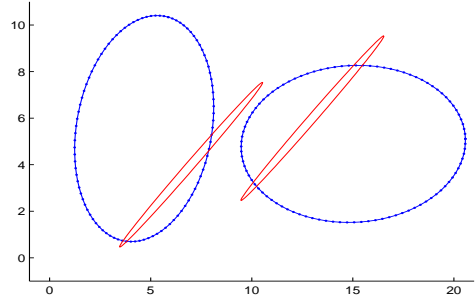}}
\subfigure[Cost as a function of algorithm iterations. $J_{e}=6.93\times10^{4}.$]{
\label{fig:T200Cost2IPA}
\includegraphics[
height=1.8in,
width=3in]
{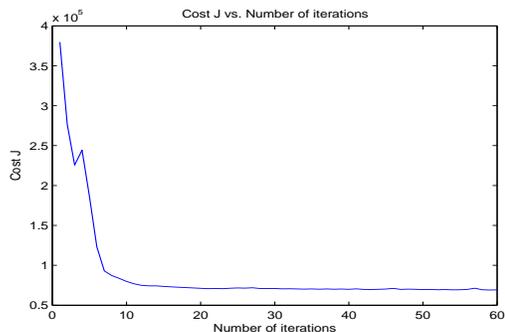}}\caption{Optimal elliptical trajectories for two agents
(without using the CSC algorithm.)}%
\label{fig:TwoAgentTraCostIPA}%
\end{figure}

We begin with a two-agent example in which we
solve \textbf{P2} by assigning elliptical trajectories using the
gradient-based approach in Section V.B (without the CSC \textbf{Algorithm 1}).
The environment setting parameters used are: $r=4$ for the sensing range of
agents; $L_{1}=20$, $L_{2}=10$, for the mission space dimensions; and $T=200$.
All sampling points $[\alpha_{i},\beta_{i}]$ are uniformly spaced within
$L_{1}\times L_{2},$ $i=1,\ldots,M$ where $M=(L_{1}+1)(L_{2}+1)=231.$ Initial
values for the uncertainty functions are $R_{i}(0)=2$ and $B=6$, $A_{i}=0.2$
for all $i=1,\ldots,M$ in (\ref{eq:multiDynR}). The results are shown in Fig.
\ref{fig:TwoAgentTraCostIPA}. Note that the initial conditions were set so as
to approximate linear trajectories (red ellipses), thus illustrating
Proposition IV.1: we can see that larger ellipses achieve a lower total
uncertainty value per unit area. Moreover, observe that the initial cost is
significantly reduced, indicating the importance of optimally selecting the
ellipse sizes, locations and orientations. The cost associated with the final
blue elliptical trajectories in this case is $J_{e}=6.93\times10^{4}$.

\begin{figure}[ptb]
\centering
\subfigure[Red and green trajectories obtained from TPBVP solution.]{
\label{fig:TPBVPcompEllipseTra}
\includegraphics[
height=1.8in,
width=3in]
{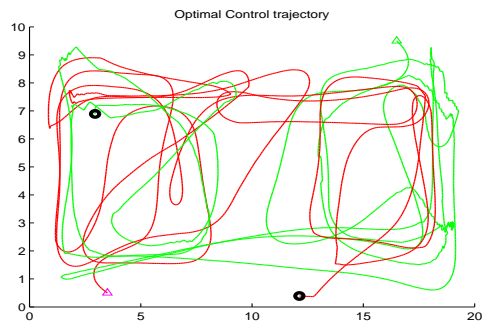}}
\subfigure[Cost as a function of algorithm iterations. $J_{\text{TPBVP}}=7.15\times10^{4}.$]{
\label{fig:TPBVPcompEllipseCost}
\includegraphics[
height=1.8in,
width=3in]
{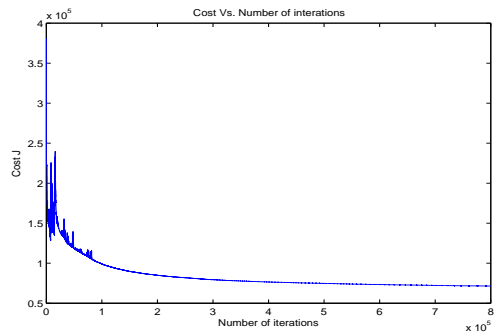}}\caption{Optimal trajectories using TPBVP solver
for two agents. Initial trajectories are red curves in Fig. \ref{fig:T200Tra}%
.}%
\label{fig:TPBVPcompEllipse}%
\end{figure}

Using the same initial trajectories as in Fig. \ref{fig:T200Tra}, we also used
a TPBVP solution algorithm for \textbf{P1}. The results are shown in Fig.
\ref{fig:TPBVPcompEllipse}. The TPBVP algorithm is computationally expensive
and time consuming (about 800,000 steps to converge). Interestingly, the
solution corresponds to a cost $J_{\text{TPBVP}}=7.15\times10^{4},$ which is
higher than that of Fig. \ref{fig:TwoAgentTraCostIPA} where solutions were
restricted to the set of elliptical trajectories. This is an indication of the
presence of locally optimal trajectories.

\begin{figure}[ptb]
\centering
\subfigure[Red ellipses: initial trajectories. Blue ellipses: optimal elliptical trajectories]{
\label{fig:T200MultiTra}
\includegraphics[
height=1.8in,
width=3in]
{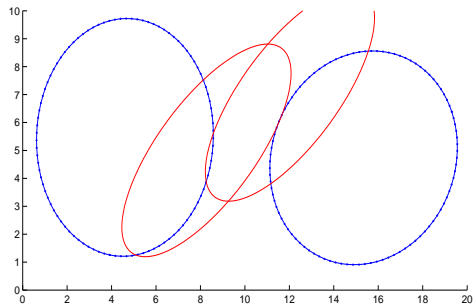}}
\subfigure[Cost as a function of algorithm iterations. $J_{\text{CSC}}^{\text{Det}} = 6.57 \times10^{4}.$]{
\label{fig:T200MultiCost}
\includegraphics[
height=1.8in,
width=3in]
{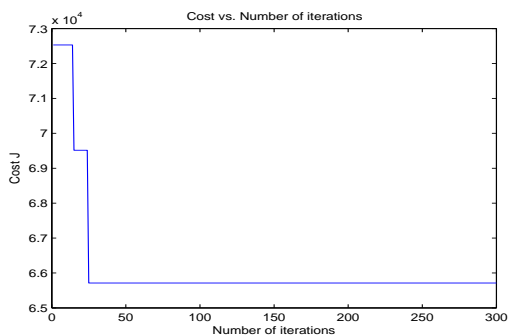}}\caption{Two agent example for the deterministic
environment setting using the CSC \textbf{Algorithm 1} for $Q=300$ trials.}%
\label{fig:T200Multi}%
\end{figure}

Next, we solve the same two-agent example with the same environment setting
using the CSC \textbf{Algorithm 1}. For simplicity, we select the ellipse
center location $[X_{n},Y_{n}]$ as the only two (out of six) multi-start
components: for a given number of comparisons $Q$, we sample the ellipse
center $[X_{n},Y_{n}]\in L_{1}\times L_{2},$ $n=1,\ldots,N,$ using a uniform
distribution while $a_{n}=5,b_{n}=2,\varphi_{n}=\frac{\pi}{4},\rho_{n}=0,$ for
$n=1,2$ are randomly assigned but initially fixed parameters during the number
of comparisons $Q$ (thus, it is still possible that there are local minima
with respect to the remaining four components $[a_{n},b_{n},\varphi_{n}%
,\rho_{n}]$, but, clearly, all six components in $\Upsilon_{n}$ can be used at
the expense of some additional computational cost.)
In Fig. \ref{fig:T200Multi}, the red elliptical trajectories on the left show
the initial ellipses and the blue trajectories represent the corresponding
resulting ellipses the CSC \textbf{Algorithm 1} converges to. Figure
\ref{fig:T200MultiCost} shows the cost vs. number of iterations of the
CSC\ algorithm. The resulting cost for $Q=300$ is $J_{\text{CSC}}^{\text{Det}%
}=6.57\times10^{4}$, where "Det" stands for a deterministic environment. It is
clear from Fig. \ref{fig:T200MultiCost} that the cost of the worst local
minimum is much higher than that of the best local minimum. Note also that the
CSC \textbf{Algorithm 1} does improve the original pure gradient-based
algorithm performance $J_{e}=6.93\times10^{4}$.

\begin{figure}[ptb]
\centering
\subfigure[Red ellipses: initial trajectories. Blue ellipses: optimal elliptical trajectories]{
\label{fig:RandomATra}
\includegraphics[
height=1.8in,
width=3in]
{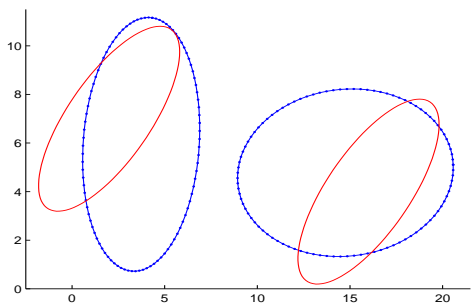}}
\subfigure[Cost as a function of algorithm iterations. $J_{\text{CSC}}^{\text{Sto}}=6.60\times10^{4}.$]{
\label{fig:RandomACost}
\includegraphics[
height=1.8in,
width=3in]
{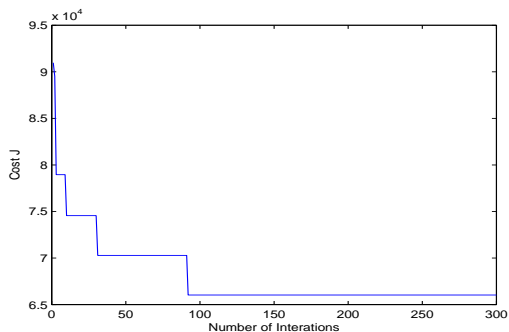}}\caption{Two-agent example for a stochastic environment
setting using the CSC \textbf{Algorithm 1} for $Q=300$ trials, where
$A_{i}\left(  \Delta t_{i}\right)  \symbol{126}U\left(  0.195,0.205\right)  $,
$\Delta t_{i}\symbol{126}0.2e^{-0.2t}$.}%
\label{fig:RandomA}%
\end{figure}

In Fig. \ref{fig:RandomA}, the values of $A_{i}$ are allowed to be
\emph{random}, thus dealing with a persistent monitoring problem in a
stochastic mission space, where we can test the robustness of the IPA approach
as discussed in \textbf{Remark 2}. In particular, each $A_{i}$ is treated as a
piecewise constant random process $\{A_{i}(t)\}$ such that $A_{i}(t)$ takes on
a fixed value sampled from a uniform distribution over $(0.195,0.205)$ for an
exponentially distributed time interval with mean $5$ before switching to a
new value. The sequence $\{M_{k}\}$ defining the number of cost comparisons
made at the $k$th iteration is set so as to grow sublinearly with
$M_{k}=\left\lceil {10\log{k}}\right\rceil ,k=2,\ldots,Q$. Note that the
system in this case is very similar to that of Fig. \ref{fig:T200Multi} where
$A_{i}=0.2$ for all $i$ without any change in the way in which $\nabla
J(\varUpsilon_{1},\ldots,\varUpsilon_{N})$ is evaluated in executing
(\ref{eq:updatetheta}). As already pointed out, this exploits a robustness
property of IPA which makes the evaluation of $\nabla J(\varUpsilon_{1}%
,\ldots,\varUpsilon_{N})$ independent of the values of $A_{i}$. All other
parameter settings are the same as in Fig. \ref{fig:T200Multi}.
In Fig. \ref{fig:RandomATra}, the red elliptical trajectories show the initial
ellipses and the blue trajectories represent the corresponding resulting
ellipses the CSC \textbf{Algorithm 1} converges to. The resulting cost for
$Q=300$ in Fig. \ref{fig:RandomACost} is $J_{\text{CSC}}^{\text{Sto}%
},=6.60\times10^{4}$, where "Sto" stands for a stochastic environment. This
cost is almost the same as $J_{\text{CSC}}^{\text{Det}}=6.57\times10^{4}$,
showing that the IPA approach is indeed robust to a stochastic environment setting.

\begin{figure}[ptb]
\centering
\subfigure[Blue ellipses: initial trajectories. Red and green trajectories: TPBVP converged trajectories. ]{
\label{fig:TPBVPAddEllipseTra}
\includegraphics[
height=1.8in,
width=3in]
{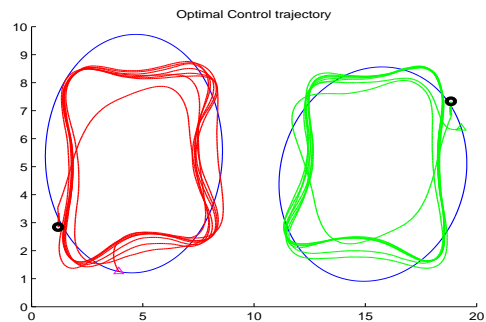}}
\subfigure[Cost vs. number of iterations. $J_{\text{TPBVP}} =
6.07\times10^{4}.$]{
\label{fig:TPBVPAddEllipseCost}
\includegraphics[
height=1.8in,
width=3in]
{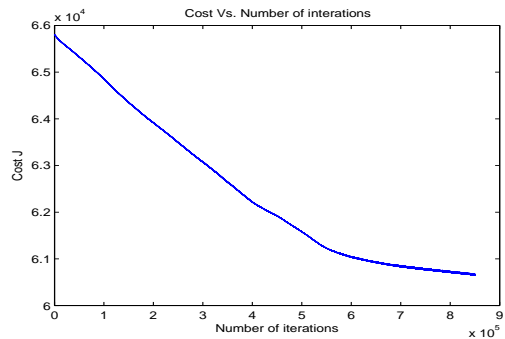}}\caption{Left plot: elliptical trajectories (blue
curve) obtained in Fig. \ref{fig:T200MultiTra} used as initial trajectories
for the TPBVP solver.}%
\label{fig:TPBVPAddEllipse}%
\end{figure}

Finally, Fig. \ref{fig:TPBVPAddEllipse} shows the TPBVP algorithm result when
using the optimal (blue) ellipses in Fig. \ref{fig:T200MultiTra} as the
initial trajectories. The trajectories the TPBVP solver converges to are shown
in red and green respectively for each agent. The corresponding cost in Fig.
\ref{fig:TPBVPAddEllipseCost} is $J_{\text{TPBVP}}=6.07\times10^{4},$ which is
an improvement compared to $J_{\text{CSC}}^{\text{Det}}=6.57\times10^{4}$
obtained for elliptical trajectories from the CSC \textbf{Algorithm 1}.
Compared to the computationally expensive TPBVP algorithm, the CSC
\textbf{Algorithm 1} using IPA is inexpensive and scalable with respect to $T$
and $N.$ Thus, a combination of the two provides the benefit of offering the
optimal elliptical trajectories obtained through the CSC \textbf{Algorithm 1}
(the first fast phase of a solution approach) as initial trajectories for the
TPBVP algorithm (the second much slower phase.) This combination is faster
than \ the original TPBVP algorithm and can also achieve a lower cost compared
to CSC \textbf{Algorithm 1}.


\section{Conclusion}

\label{sec:concl} We have shown that an optimal control solution to the 1D
persistent monitoring problem does not easily extend to the 2D case. In
particular, we have proved that elliptical trajectories outperform linear ones
in a 2D mission space. Therefore, we have sought to solve a parametric
optimization problem to determine optimal elliptical trajectories. Numerical
examples indicate that this scalable approach (which can be used on line)
provides solutions that approximate those obtained through a computationally
intensive TPBVP solver. Moreover, since the solutions obtained are generally
locally optimal, we have incorporated a stochastic comparison algorithm for
deriving globally optimal elliptical trajectories. Ongoing work aims at
alternative approaches for near-optimal solutions and at distributed implementations.

\bibliographystyle{IEEEtran}
\bibliography{Papers}

\end{document}